\documentclass[12pt]{article}
\usepackage[cp1251]{inputenc}

\usepackage[english]{babel}
\usepackage{amsmath,amsfonts,amscd,amssymb,latexsym}

\makeindex
\begin{document}
\begin{center}
{\Large  Structures of Malcev Bialgebras on a simple non-Lie Malcev
algebra.}
\end{center}

\vspace{7mm}

\begin{center}
{\bf M.\,E.\,Goncharov}
\end{center}

\footnotetext{The author was supported by Lavrent'ev Young
Scientists Competition (No 43 on 04.02.2010), ADTP ``Development of
the Scientific Potential of Higher School'' of the Russian Federal
Agency for Education (Grant 2.1.1.419), Russian Foundation for basic
Research (Grant 09-01-00157-A), the State Support Programs for
the Leading Scientific Schools and the Young Scientists of the
Russian Federation (Grands Nsh-3669.2010.1 and MD-2438.2009.1), and
the Federal Targeting Programs (contracts 02.740.11.0429,
02.740.11.5191), the Integration Grant of the Siberian Division of
the Russian Academy of Sciences (No. 97). }

\renewcommand{\abstractname}{Abstract}
\begin{abstract}
In this work, we consider Malcev bialgebras. We describe all structures
of a Malcev bialgebra on a simple non-Lie Malcev algebra.

\end{abstract}
 \vspace{6pt}

Lie bialgebras are Lie algebras and Lie coalgebras at the same time,
such that comultiplication is a 1-cocycle. These bialgebras
were introduced by Drinfeld [1] in studying the solutions to the
classical Yang-Baxter equation. In [2, 3], the
definition of a bialgebra in the sense of Drinfeld (D-bialgebra)
related with any variety of algebras was stated.
In particular, the
associative and Jordan D-bialgebras were introduced, and an
associative analogue of the Yang-Baxter equation was considered as
well as the associative D-bialgebras related with the solutions to
this equation.
In the same papers, the associative algebras
that admit a nontrivial structure of a D-bialgebra with
cocommutative comultiplication on the center were also
described.
The comultiplication
in an associative D-bialgebra is a derivation of an initial
algebra into its tensor square considered as a bimodule over the
initial algebra.
These bialgebras were introduced in [4] and studied
in [5]. The paper [5] is devoted to some properties  of solutions to an
associative analogue of the Yang-Baxter equation and the balanced
bialgebras (i.e., D-bialgebras).
The associative classical
Yang-Baxter equations with parameters were considered in [6]. A
class of Jordan D-bialgebras related to the ``Jordan analogue'' of the
Yang-Baxter equation was introduced in [7], where it was proved that
every finite-dimensional Jordan D-bialgebra, semisimple as an
algebra, belongs to this class.

A so-called Manin triple may be associated with every Lie,
associative, or Jordan bialgebra. In [8], the Manin triples for the
associative algebras served as a tool for the study of the solutions
to the Yang-Baxter equation.

Alternative D-bialgebras and their connection with the alternative
Yang-Baxter equation were under study in [9]. In particular, the
alternative D-bialgebra structures on Cayley-Dickson matrix algebra
were described. Some connection of Jordan D-bialgebras with Lie
bialgebras was revealed in [2]. It was shown in particular that
under some natural restrictions a Jordan algebra $J$ admits a
nontrivial structure of a Jordan D-bialgebra if the Lie algebra
$L(J)$ obtained by the Kantor-Koecher-Tits (KKT) construction from
$J$ admits the structure of a Lie bialgebra. Given an associative
D-bialgebra $(A, \Delta)$ and its adjoint Jordan D-bialgebra
$(A^{(+)},\Delta^{(+)})$, $L(A^{(+)})$ may be equipped with the
structure of a Lie bialgebra which is connected in some sense with
$(A^{(+)},\Delta^{(+)})$ [2]. In the present paper, we prove an
analogue to this result in the case when $A$ is a Cayley-Dickson
matrix algebra and $(A,\Delta)$ is an alternative D-bialgebra. At
the same time, we construct an example of alternative D-bialgebra
$(A, \Delta)$ for which the structure of an adjoint Jordan
D-bialgebra on $(A^{(+)},\Delta^{(+)})$ cannot be extended to the
structure of a Lie bialgebra on $L(A^{(+)})$.

Functional solutions to the classical Yang-Baxter equation on simple
Lie algebras were constructed in \cite{BD}. In \cite{stolin}, using
the ideas from \cite{BD}, it was obtained an explicit description of
Lie bialgebra structures on simple complex Lie algebras.

Malcev algebras were introduced by A.I. Malcev \cite{M55} as tangent
algebras for local analytic Moufang loops.  The class of Malcev
algebras generalizes the class of Lie algebras and has a well
developed theory \cite{KuzSh}.

An important example of a non-Lie Malcev algebra is the vector space
of zero trace elements of a Caley-Dickson algebra with the commutator
bracket multiplication \cite{S62,K68}. In \cite{versh} some
properties of Malcev bialgebras were studied. In particular, there were
found  conditions for a Malcev algebra with a comultiplication to
be a Malcev bialgebra.

In this work, we consider an analogue of the classical Yang-Baxter
equation on Malcev algebras. In particular, it is shown that any
solution to this equation induces a structure of a Malcev bialgebra.
Also, we describe all structures of a Malcev bialgebra on the simple
non-Lie complex Malcev algebra.

In order to perform vast volume of routine computations we used the
Groups, Algorithms, Programming  System (GAP).

\begin{center}
{\bf \S 1. Definitions and Preliminaries}
\end{center}
Given vector spaces $V$ and $U$ over a field $F$, denote by
$V\otimes U$ its tensor product over $F$. Define the linear mapping
$\tau$ on $V$ by  $\tau(\sum\limits_ia_i\otimes
b_i)=\sum\limits_ib_i\otimes a_i$. Define the linear mapping $\xi$
on $V\otimes V\otimes V$ by  $\xi(\sum\limits_i a_i\otimes
b_i\otimes c_i)= \sum\limits_i b_i\otimes c_i\otimes a_i$. Denote by
$V^*$ the dual space of $V$. Given $f\in V$ and
$\upsilon \in V$, the symbol $\langle f, \upsilon\rangle$ denotes
the linear functional $f$ evaluated at $\upsilon$ (i.e., $\langle
f,a\rangle=f(a)$).

{\bf Definition.} A pair $(A, \Delta)$, where $A$ is a vector space
over $F$ and $\Delta : A \rightarrow A\otimes A$ is a linear
mapping, is called a \emph{coalgebra}, while $\Delta$ is a
\emph{comultiplication}.

Given $a\in A$, put $\Delta(a) = \sum a_{(1)}\otimes a_{(2)}$.

Define some multiplication on $A^*$ by
$$
\langle fg,a\rangle=\sum\limits_a\langle f,a_{(1)}\rangle\langle
g,a_{(2)}\rangle,
$$
where $f,g\in A^*$, $a\in A$ and $\Delta(a)=\sum a_{(1)}\otimes
a_{(2)}$.  The algebra obtained is  \emph{the dual algebra} of the
coalgebra $(A, \Delta)$.

The dual algebra $A^*$ of $(A,\Delta)$ gives rise to the following
bimodule actions on $A$:
$$
f\rightharpoonup a=\sum a_{(1)}\langle f,a_{(2)}\rangle\text{ and
}a\leftharpoonup f=\sum\langle f,a_{(1)}\rangle a_{(2)},
$$
where $a\in A,\ f\in A^*$ and $\Delta(a)=\sum a_{(1)}\otimes
a_{(2)}$.

The following definition  of a coalgebra related to some variety of
algebras was given in \cite{ANQ}.

{\bf Definition.} Let $\mathcal{M}$ be an arbitrary variety of
algebras. The pair $(A, \Delta)$ is called a $\mathcal{M}$-coalgebra
if $A^*$ belongs to $M$.

Let $A$ be an arbitrary algebra with a comultiplication $\Delta$, and
let $A^*$ be the dual algebra for $(A, \Delta)$. Then $A$ induces
the bimodule action on $A^*$ by the formulas
$$
\langle f\leftharpoondown a,b\rangle=\langle f,ab\rangle\text{ and }
\langle b\rightharpoondown f,a\rangle=\langle f,ab\rangle,
$$
where $a,b\in A,\ f\in A^*$.

Consider the space $D(A) = A \oplus A$ and equip it with the
multiplication by putting
$$
(a+f)(b+g)=(ab+f\rightharpoonup b+a\leftharpoonup
g)+(fg+f\leftharpoondown b+a\rightharpoondown g).
$$
Then $D(A)$ is an ordinary algebra over $F$, $A$ and $A^*$ are some
subalgebras in $D(A)$. It is called \emph{the Drinfeld double}.

Let $Q$ be a bilinear form on $D(A)$ defined by
$$
Q(a+f,b+g)=\langle g,a\rangle+\langle f,b\rangle
$$
for all $a,b\in A$ and $f,g\in A^*$. It is easy to check that $Q$ is
a nondegenerate symmetric associative form, that is
$Q(xy,z)=Q(x,yz)$.

Let us recall the definition of a Lie bialgebra. Let $L$ be a Lie
algebra with a comultiplication $\Delta$. The pair $(A,\Delta)$ is
called a Lie bialgebra if and only if $(L,\Delta)$ is a Lie
coalgebra and $\Delta$ is a 1-cocycle, i.e., it satisfies
$$
\Delta([a,b])=\sum([a_{(1)},b]\otimes
a_{(2)}+a_{(1)}\otimes[a_{(2)},b])+\sum([a,b_{(1)}]\otimes
b_{(2)}+b_{(1)}\otimes[a,b_{(2)}])
$$
for all $a,b\in L$.

In [1], it was proved that the pair $(L,\Delta)$ is a Lie bialgebra
if and only if its Drinfeld double $D(L)$ is a Lie algebra.
This observation inspired the following definition~[2].

{\bf Definition.} Let $\mathcal{M}$ be an arbitrary variety of
algebras and let $A$ be an algebra from $\mathcal{M}$ with a
comultiplication $\Delta$. The pair $(A,\Delta)$ is called \emph{an
$\mathcal{M}$-bialgebra in the sense of Drinfeld\/} if its Drinfeld
double $D(A)$ belongs to $\mathcal{M}$.

Note that this definition corresponds with the definition of
coalgebra given in \cite{ANQ}.

There is an important type of Lie bialgebras
called coboundary bialgebras. Namely, let $L$ be a Lie
algebra and $R=\sum\limits_ia_i\otimes b_i$ from $(id-\tau)(L\otimes
L)$, that is, $\tau(R)=-R$. Define a comultiplication $\Delta_R$ on
$L$ by
$$
\Delta_R(a)=\sum\limits_i[a_i,a]\otimes b_i-a_i\otimes[a,b_i]
$$
for all $a\in L$. It is easy to see that $\Delta_R$ is a 1-cocycle.
In [11] it was proved that $(L,\Delta)$ is a Lie coalgebra if and
only if the element
$$
C_L(R)=[R_{12},R_{13}]+[R_{12},R_{23}]+[R_{13},R_{23}]
$$
is $L$-invariant. Here
$[R_{12},R_{13}]=\sum\limits_{ij}[a_i,a_j]\otimes b_i\otimes b_j$,
$[R_{12},R_{23}]=\sum\limits_{ij}b_i\otimes[a_i,a_j]\otimes b_j$,
 and $[R_{13},R_{23}]=\sum\limits_{ij} a_i\otimes a_j\otimes [b_i,b_j]$.
In particular, if
\begin{equation}\label{lieYB}
[R_{12},R_{13}]+[R_{12},R_{23}]+[R_{13},R_{23}]=0,
\end{equation}
then the pair $(L,\Delta_R)$ is a Lie bialgebra. In this case, we say that
$(L,\Delta_R)$ is a triangular Lie bialgebra. The equation
\eqref{lieYB} is called \emph{the classical Yang-Baxter equation}.


Let $B$ be an arbitrary algebra and $r=\sum\limits_ia_i\otimes
b_i\in B\otimes B$. Then the equation
\begin{equation}\label{YB}
  C_B(r)=r_{12}r_{13}+r_{13}r_{23}-r_{23}r_{12}=0
\end{equation}
is called the classical Yang-Baxter equation on $B$. Here the
subscripts specify the way of embedding $B\otimes B$ into $B\otimes
B\otimes B$, that is, $r_{12}=\sum\limits_{i} a_i\otimes b_i\otimes
1$, $r_{13}=\sum_i a_i\otimes 1\otimes b_i$, $r_{23}=\sum_i 1\otimes
a_i\otimes b_i$. Note that $C_B(r)$ is well defined even if $B$ is
non-unital. This equation for different varieties of algebras were
considered in \cite{Zhelyabin98, Aquiar, Zhelyabin, Polishchuk,
Gme}. Usually, antisymmetric solutions ($\tau(r)=-r$) to the
equation \eqref{YB} are considered.

An element $r=\sum\limits_{i}a_i\otimes b_i\in B\otimes B$ induces a
comultiplication $\Delta_r$ on $B$:
$$
\Delta_r(a)=\sum\limits a_ia\otimes b_i-a_i\otimes ab_i
$$
for all $a\in B$.

We will need the following

{\bf Lemma 1.} Let $B$ be an arbitrary anticommutative
finite-dimensional algebra over a field $F$, $a_1,\ldots,a_n$ be a
basis of $B$, and $\gamma_{ij}^k$ be the structure constants of $B$
with respect to $a_1,\ldots,a_n$, i.e.,  $a_ia_j=\sum\limits_k
\gamma_{ij}^ka_k$. Suppose that for an element $r=\sum\limits_{ij}
\alpha_{ij} a_i\otimes a_j \in B\otimes B$ and the dual algebra
$B^*$ of coalgebra $(B,-\Delta_r)$ the mapping $\phi: B^*\rightarrow
B$ defined by $\phi(f)=\sum\limits_{ij}\alpha_{ij}f(a_j)a_i$ is a
homomorphism of algebras.  Then the following equations hold:
\begin{equation}\label{eq2}
(\Lambda^{\top}\Gamma_k\Lambda)_{sn}+(\Lambda\Gamma_s\Lambda)_{kn}+
(\Lambda\Gamma_n\Lambda^{\top})_{ks}=0,
\end{equation}
\begin{equation}\label{nd}
\det(\Lambda)(\sum_l2(\Lambda\Gamma_l)_{kl}+(\Lambda^{\top}\Gamma_k)_{ll})=0,
\end{equation}
where $\Lambda=(\alpha_{ij})_{i,j=1\ldots n}$,
$\Gamma_k=(\gamma_{ij}^k)_{i,j=1\ldots n}$, and $\Lambda^{\top}$ is the
matrix transpose to $\Lambda$.

{\bf Proof.} Let $b_k=\sum\limits_i \alpha_{ki}a_i$. Then
$r=\sum\limits_i a_i\otimes b_i$. Since the mapping $\phi$ is a
homomorphism, for all $f,g\in A^*$ we have
$$
\sum_{i,j}f(a_i)g(a_j)b_ib_j=\sum_ifg(a_i)b_i=\sum\limits_i\langle
f\otimes g,-\Delta_r(a_i)\rangle
b_i=
$$
$$
=-(\sum\limits_{ij}f(a_ja_i)g(b_j)b_i-\sum\limits_{ij}f(a_j)g(a_ib_j)b_i).
$$
Thus $r$ satisfies  the following condition:
\begin{equation}\label{YB1}
\sum\limits_{ij}a_ia_j\otimes b_i\otimes b_j-a_i\otimes
a_jb_i\otimes b_j+a_i\otimes a_j\otimes b_ib_j=0,
\end{equation}
that is, $r$ is a solution to the classical Yang-Baxter equation.
Rewriting the elements $b_i$ in terms of the elements $a_i$, one can
get
$$
\sum\limits_{i,j,n,s,k}\gamma_{ij}^k\alpha_{is}\alpha_{jn}a_k\otimes
a_s\otimes a_n -\gamma_{js}^k\alpha_{is}\alpha_{jn} a_i\otimes
a_k\otimes a_n+\gamma_{ns}^k\alpha_{in}\alpha_{js}a_i\otimes
a_j\otimes a_k=0.$$

Changing corresponding indices in the second and third summands and
using the skew-symmetry of $\Gamma_i$ we conclude

$$
\sum\limits_{i,j,n,s,k}(\gamma_{ij}^k\alpha_{jn})\alpha_{is}a_k\otimes
a_s\otimes a_n +(\alpha_{ki}\gamma_{ij}^s)\alpha_{jn} a_k\otimes
a_s\otimes a_n+(\alpha_{ki}\gamma_{ij}^n)\alpha_{sj}a_k\otimes
a_s\otimes a_n=0.$$

Therefore, for all $k,s,n$
$$
\sum\limits_{ij}
(\gamma_{ij}^k\alpha_{jn})\alpha_{is}+(\alpha_{ki}\gamma_{ij}^s)\alpha_{jn}+(\alpha_{ki}\gamma_{ij}^n)\alpha_{sj}=0$$

Changing indices $i$ and $j$ in the first summand we obtain

$$
\sum\limits_{ij}
(\gamma_{ji}^k\alpha_{in})\alpha_{js}+(\alpha_{ki}\gamma_{ij}^s)\alpha_{jn}+(\alpha_{ki}\gamma_{ij}^n)\alpha_{sj}=0$$

Hence,
\begin{equation}\label{eq1} \sum\limits_j
(\Gamma_k\Lambda)_{jn}\alpha_{js}+(\Lambda\Gamma_s)_{kj}\alpha_{jn}+(\Lambda\Gamma_n)_{kj}\alpha_{sj}=0.
\end{equation}
One can rewrite this equation in the following form:
$$
(\Lambda^{\top}\Gamma_k\Lambda)_{sn}+(\Lambda\Gamma_s\Lambda)_{kn}+(\Lambda\Gamma_n\Lambda^{\top})_{ks}=0.
$$
This  proves the first statement \eqref{eq2}.

Let $\Lambda^*=(\alpha_{ij}^*)$ be the cofactor matrix of $\Lambda$.
Then  $\Lambda\Lambda^{*\top}=\det(\Lambda)E,$ where  $E$ is the
identity matrix. Multiplying the equation \eqref{eq1} on
$\alpha^*_{ln}$ and summing up over $n$ we get

$$
\sum\limits_{jn}(\Gamma_k\Lambda)_{jn}\alpha_{js}\alpha^{*}_{ln}+(\Lambda\Gamma_s)_{kj}\alpha_{jn}\alpha^*_{ln}+(\Lambda\Gamma_n)_{kj}\alpha_{sj}\alpha^*_{ln}=0.
$$
Therefore,
$$
\sum\limits_{jn}((\Gamma_k\Lambda)_{jn}\alpha_{js}\alpha^{*}_{ln}+(\Lambda\Gamma_n)_{kj}\alpha_{sj}\alpha^*_{ln})+(\Lambda\Gamma_s)_{kl}\det(\Lambda)=0.
$$
The last equation can be rewritten in the form
$$
\sum\limits_{n}((\Lambda^{\top}\Gamma_k\Lambda)_{sn}\alpha^*_{ln}+(\Lambda\Gamma_n\Lambda^{\top})_{ks}\alpha^*_{ln})+(\Lambda\Gamma_s)_{kl}\det(\Lambda)=0.
$$
Thus
$$
\sum\limits_{n}(\Lambda\Gamma_n\Lambda^{\top})_{ks}\alpha^*_{ln}+(\Lambda^{\top}\Gamma_k\Lambda\Lambda^{*\top})_{sl}+(\Lambda\Gamma_s)_{kl}\det(\Lambda)=0.
$$
Hence,
$$
\sum\limits_{n}(\Lambda\Gamma_n\Lambda^{\top})_{ks}\alpha^*_{ln}+(\Lambda^{\top}\Gamma_k)_{sl}\det(\Lambda)+(\Lambda\Gamma_s)_{kl}\det(\Lambda)=0.
$$
For $s=l$ we have
$$
\sum\limits_{n}(\Lambda\Gamma_n\Lambda^{\top})_{kl}\alpha^*_{ln}+(\Lambda^{\tau}\Gamma_k)_{ll}\det(\Lambda)+(\Lambda\Gamma_l)_{kl}\det(\Lambda)=0.
$$
Summing up over $l$ we obtain
$$
\sum\limits_l\sum\limits_{n}(\Lambda\Gamma_n\Lambda
A^{\top})_{kl}\alpha^*_{ln}+\sum\limits_l(\Lambda^{\top}\Gamma_k)_{ll}\det(\Lambda)+\sum\limits_l(\Lambda\Gamma_l)_{kl}\det(\Lambda)=0.
$$
Then
$$
\sum\limits_{n}(\Lambda\Gamma_n\Lambda
^{\top}\Lambda^*)_{kn}+\sum\limits_l(\Lambda^{\top}\Gamma_k)_{ll}\det(\Lambda)+\sum\limits_l(\Lambda\Gamma_l)_{kl}\det(\Lambda)=0.
$$
Therefore,
$$
\sum\limits_{n}(\Lambda\Gamma_n\det(\Lambda))_{kn}+\sum\limits_l(\Lambda^{\top}\Gamma_k)_{ll}\det(\Lambda)+\sum\limits_l(\Lambda\Gamma_l)_{kl}\det(\Lambda)=0.
$$
Finally, changing indexes $n$ on $l$ in the first sum, we obtain
$$
\det(\Lambda)(\sum_l2(\Lambda\Gamma_l)_{kl}+(\Lambda^{\top}\Gamma_k)_{ll})=0.
$$\hspace*{\fill} $\Box$

The following lemma was originally proved in \cite{Zhelyabin}. Here
we state its proof in order to complete the exposition.

 {\bf Lemma 2.} Let $B$ be a finite-dimensional simple algebra over a field $F$
 with a nontrivial comultiplication $\Delta$,
  $D(B)=B\oplus B^*$ be the Drinfeld double of the coalgebra  $(B,\Delta)$.
Suppose that $U$ is a nonzero ideal in  $D(B)$ and
$$
V=\{a\in B|a+f\in U\ \text{for some}\ f\in B^*\}.
$$
Then the dimensions of  $B$ and $U$ are equal, the pair $(V,
\Delta)$ is a subbialgebra of $(B,\Delta)$, and
$V^{\perp}U=UV^{\perp}=0$, where $V^{\perp}$ is the orthogonal
complement of  $V$ in $B^*$.

\textsc{Proof.} Let $N$ be a $B$-subbimodule in $B^*$ and
$N^{\perp}$ be the orthogonal complement of $N$ in $B$ with respect
to the form $Q$. Then  $N^{\perp}$ is an ideal in $J$. Since $B$ is
a simple algebra, then either $N^{\perp}=0$ or $N^{\perp}=B^*$.
Therefore $B^*$ is an irreducible $B$-bimodule.

Consider the vector space $W=\{ f\in B^*|f+a\in U\ \text{ for some}\
a\in B\}$. Since $D(B)=B\oplus B^*$, then $W\neq0$. Take $f\in W$
and $b\in B$. Then for some $a$ from $B$ we have $f+a\in U$. Since
$(f+a)b=f\rightharpoonup b+ab+f\leftharpoondown b\in U$, then
$f\leftharpoondown b\in W$. Similarly,  $b\rightharpoondown f\in W$.
Therefore $W$ is a $B$-subbimodule in $B^*$ and $W=B^*$. Hence,
$dim_FU\geq dim_FB^*$. Note that $U\cap B=0$ since $B$ is a simple
algebra. Therefore, $dim_F(U+B)=dim_FD(B)$ and $dim_F
U=dim_FB^*=dim_FB$.

Now, let us show that the pair $(V,\Delta)$ is a subbialgebra in
$(B,\Delta)$. For this, take $a\in V$ and $g\in B^*$. Then, for some
$f\in B^*$ we have $a+f\in U$. Hence
$$(a+f)g=fg+a \rightharpoondown g + a \leftharpoonup g \in U.$$
Therefore $a\leftharpoonup g\in V$. Similarly, $g\rightharpoonup a\in
V$ and $(V,\Delta)$ is a subcoalgebra in  $(B,\Delta)$. Consider
$b\in V$, then
$$(a+f)b=f\leftharpoondown b+f\rightharpoonup b+ab\in
U.$$
Hence $V$ is a subbialgebra in $B$ and the pair $(V,\Delta)$ is
a subbialgebra in $(B,\Delta)$.

Let $b,n$ be arbitrary elements from $B$ and $V^{\perp}$
respectively. Then $$(a+f)b=f\leftharpoondown b + f \rightharpoonup
b+ab\in U.$$ Therefore $f\rightharpoonup b+ab\in V$. This implies
 that $Q(nf+n\leftharpoondown
a,b)=Q(n,f\rightharpoonup b+ab)=0$. Thus, we obtain
$nf+n\leftharpoondown a=0$. But any element $u\in U$ can be
represented in the form $f+a$, where $f\in B^*$ and $a\in V$. Hence,
$nu=n(a+f)=nf+n\leftharpoondown a+n\rightharpoonup a=0$ and this
proves that $V^{\perp}U=0$. Similarly, $UV^{\perp}=0$.
\hspace*{\fill} $\Box$

\begin{center}
{\bf \S 2. Coboundary and triangular Malcev bialgebras.}
\end{center}

An anticommutative algebra is called a Malcev algebra if for all
$x,y,z\in M$ the following equation holds:
\begin{equation}\label{mal1}
J(x,y,xz)=J(x,y,z)x,
\end{equation}
where  $J(x,y,z)=(xy)z+(yz)x+(zx)y$ is the jacobian of elements
$x,y,z$.

Note that if the characteristic of the field $F$ is different from
2 then  \eqref{mal1} is equivalent to the following equation:
\begin{equation}\label{mal}
((xy)z)t+((yz)t)x+((zt)x)y+((tx)y)z=(xz)(yt).
\end{equation}

For Malcev bialgebras, one can also consider the class of coboundary
and triangular bialgebras.

Let $M$ be a Malcev algebra, $r\in(id-\tau)(M\otimes M)$. Then $r$
induces a comultiplication $\Delta_r$ on $M$ defined by
$$
\Delta_r(a)=[r,a]=\sum\limits_ia_ia\otimes b_i-a_i\otimes ab_i
$$
for all $a\in M$. In this work, we find necessary and sufficient
conditions for the pair $(M,\Delta_r)$ to be a Malcev bialgebra.

In \cite{versh}, the following statement was proved.

{\bf Theorem 1.(Vershinin).} A Malcev algebra $M$ with a
multiplication $\mu$ and a comultiplication $\Delta$ is a Malcev
bialgebra if and only if the dual algebra $M^*$ is a Malcev algebra
and the comultiplication satisfies
\begin{enumerate}
\item$\Delta((ab)c)+\Delta(bc)(a\otimes 1)+(b\otimes
1)\Delta(ac)=$ $\sum a_{(1)}(bc)\otimes a_{(2)}+\sum a_{(1)}c\otimes
a_{(2)}b+\sum a_{(1)}\otimes (a_{(2)}b)c+$ $\sum\limits
ab_{(1)}\otimes b_{(2)}c-\sum b_{(1)}\otimes
b_{(2)}(ac)+\sum(ab_{(1)})c\otimes b_{(2)}$ $+\sum
a(bc_{(1)})\otimes c_{(2)}-\sum c_{(1)}\otimes (c_{(2)}a)b ,$

\item$(1\otimes \Delta)\Delta(ab)= (1\otimes 1\otimes a)((1\otimes
\Delta)\Delta(b))+(\Delta\otimes 1)((a\otimes 1)\Delta(b))-(1\otimes
\tau)((1\otimes 1\otimes a)(\Delta\otimes 1)\Delta(b))+(1\otimes
\tau)(((\Delta\otimes 1)\Delta(b))(a\otimes 1\otimes
1))+((\Delta\otimes 1)\Delta(a))(b\otimes 1\otimes
1)+((\Delta\otimes 1)\Delta(a))(1\otimes 1\otimes b) -(1\otimes
\tau)((\Delta\otimes 1)(\Delta(a)(b\otimes 1)))-(1\otimes b\otimes
1)((1\otimes \Delta)\Delta(a))+(1\otimes \tau)(1\otimes 1\otimes
\mu)((1\otimes \tau\otimes 1)(\Delta(b)\otimes\Delta(a)))+(1\otimes
1\otimes \mu)((1\otimes \tau\otimes 1)(\Delta(a)\otimes
\Delta(b))).$
\end{enumerate}

In this work, we prove

{\bf Theorem 2.} Let $M$ be a Malcev algebra over a field
characteristic not equal 2, $r\in (id-\tau)(M\otimes M)$. The pair
$(M,\Delta_r)$ is a Malcev bialgebra if and only if for all
$a,b\in M$
$$
(C_M(r)(1\otimes b\otimes 1))(1\otimes a\otimes 1)-C_M(r)(ab\otimes
1\otimes 1)-(C_M(r)(1\otimes 1\otimes a)),(1\otimes 1\otimes b)=
$$
\begin{equation}\label{UM}
=C_M(r)(b\otimes 1\otimes a)-C_M(r)(a\otimes b\otimes 1)
\end{equation}
or
$$
C_M(r)(1\otimes J_{b,a}\otimes 1-1\otimes 1\otimes
J_{a,b})=[C_M(r),ab]+[C_M(r),b](1\otimes 1\otimes
a)-[C_M(r),a](1\otimes b\otimes 1),
$$
where the operator  $J_{a,b}$ is defined by $cJ_{a,b}=J(c,a,b)$, and
by $[C_M(r),a]$ we denote the action of $M$ on $M\otimes M\otimes M$
defined by $[x\otimes y\otimes z,a]=xa\otimes y\otimes z+x\otimes
ya\otimes z+x\otimes y\otimes za$.

\textsc{Proof.} Let us prove the necessary condition. Since the pair
$(M,\Delta_r)$ is a Malcev bialgebra then the second equation of the
theorem 1 holds.
We have:\\
 $ (1\otimes \Delta)\Delta_r(ab)=$
$$\sum\limits_{ij}a_i(ab)\otimes a_jb_i\otimes
b_j-\sum\limits_{ij}a_i(ab)\otimes a_j\otimes
b_ib_j-\sum\limits_{ij}a_i\otimes a_j((ab)b_i)\otimes b_j
+\sum\limits_{ij}a_i\otimes a_j\otimes((ab)b_i)b_j.\\
$$
$(1\otimes 1\otimes a)((1\otimes \Delta)\Delta(b))=$
$$
\sum\limits_{ij}a_ib\otimes a_jb_i\otimes
ab_j-\sum\limits_{ij}a_ib\otimes a_j\otimes
a(b_ib_j)-\sum\limits_{ij}a_i\otimes a_j(bb_i)\otimes
ab_j+\sum\limits_{ij}a_i\otimes a_j\otimes a((bb_i)b_j).\\
$$
$(\Delta\otimes 1)((a\otimes 1)\Delta(b))=$
$$
\sum\limits_{ij}a_j(a(a_ib))\otimes b_j\otimes b_i-\sum\limits_{ij}
a_j\otimes (a(a_ib))b_j\otimes b_i-\sum\limits_{ij}a_j(aa_i)\otimes
b_j\otimes bb_i+a_j\otimes(aa_i)b_j\otimes bb_i.\\$$ $ -(1\otimes
\tau)((1\otimes 1\otimes a)(\Delta\otimes 1)\Delta(b))=$
$$
-\sum\limits_{ij}a_j(a_ib)\otimes ab_i\otimes
b_j+\sum\limits_{ij}a_ja_i\otimes a(bb_i)\otimes
b_j+\sum\limits_{ij}a_j\otimes ab_i\otimes
(a_ib)b_j-\sum\limits_{ij}a_j\otimes a(bb_i)\otimes a_ib_j.\\$$ $
(1\otimes \tau)(((\Delta\otimes 1)\Delta(b))(a\otimes 1\otimes 1))=$
$$\sum\limits_{ij}(a_j(a_ib))a\otimes b_i \otimes b_j
-\sum\limits_{ij}a_ja\otimes b_i\otimes
(a_ib)b_j-\sum\limits_{ij}(a_ja_i)a\otimes bb_i\otimes
b_j+\sum\limits_{ij}a_ja\otimes a_ib_j\otimes bb_i.$$
 $
((\Delta\otimes 1)\Delta(a))(b\otimes 1\otimes 1)=$
$$
\sum\limits_{ij}(a_j(a_ia))b\otimes b_j\otimes
b_i-\sum\limits_{ij}a_jb\otimes(a_ia)b_j\otimes
b_i-\sum\limits_{ij}(a_ja_i)b\otimes b_j\otimes
ab_i+\sum\limits_{ij}a_jb\otimes a_ib_j\otimes ab_i.\\$$ $
((\Delta\otimes 1)\Delta(a))(1\otimes 1\otimes b)=$
$$
\sum\limits_{ij} a_j(a_ia)\otimes b_j\otimes
b_ib-\sum\limits_{ij}a_j\otimes(a_ia_j)b_j\otimes
b_ib-\sum\limits_{ij}a_ja_i\otimes b_j\otimes
(ab_i)b+\sum\limits_{ij}a_j\otimes a_ib_j\otimes (ab_i)b.\\
$$
 $ -(1\otimes
\tau)((\Delta\otimes 1)(\Delta(a)(b\otimes 1)))=$
$$
-\sum\limits_{ij} a_j((a_ia)b)\otimes b_i\otimes
b_j+\sum\limits_{ij}a_j\otimes b_i\otimes
((a_ia)b)b_j+\sum\limits_{ij}a_j(a_ib)\otimes ab_i\otimes
b_j-\sum\limits_{ij}a_j\otimes ab_i\otimes(a_ib)b_j.\\
$$
$ -(1\otimes b\otimes 1)((1\otimes \Delta)\Delta(a))=$
$$-\sum\limits_{ij}a_ia\otimes b(a_jb_i)\otimes
b_j+\sum\limits_{ij}a_ia\otimes ba_j\otimes
b_ib_j+\sum\limits_{ij}a_i\otimes b(a_j(ab_i))\otimes
b_j-\sum\limits_{ij}a_i\otimes ba_j\otimes (ab_i)b_j.\\$$
 $(1
\otimes \tau)(1\otimes 1\otimes \mu)((1\otimes\tau\otimes
1)(\Delta(b)\otimes \Delta(a)))=$
$$
\sum\limits_{ij}a_ib\otimes b_ib_j\otimes a_ja-
-\sum\limits_{ij}a_ib\otimes b_i(ab_j)\otimes
a_j-\sum\limits_{ij}a_i\otimes(bb_i)b_j\otimes
a_ja+\sum\limits_{ij}a_i\otimes (bb_i)(ab_j)\otimes a_j.\\$$ $
(1\otimes1\otimes\mu)((1\otimes\tau\otimes 1)(\Delta(a)\Delta(b)))=$
$$\sum\limits_{ij}a_ia\otimes a_jb\otimes
b_ib_j-\sum\limits_{ij}a_i\otimes a_jb\otimes
(ab_i)b_j-\sum\limits_{ij}a_ia\otimes a_j\otimes
b_i(bb_j)+\sum\limits_{ij}a_i\otimes a_j\otimes (ab_i)(bb_j).\\$$

Inserting the expressions obtained into second equality of the
theorem 1, using \eqref{mal} and $\tau(r)=-r$, we conclude:

$$
\sum\limits_{ij} a_i(ab)\otimes a_jb_i\otimes b_j
-\sum\limits_{ij}a_i(ab)\otimes a_j\otimes
b_ib_j-\sum\limits_{ij}a_i\otimes((b_ib_j)b)a\otimes
a_j-\sum\limits_{ij}a_i\otimes a_j\otimes ((b_ib_j)a)b=
$$
$$
=-\sum\limits_{ij}a_ib\otimes a_j\otimes
a(b_ib_j)+\sum\limits_{ij}a_ib\otimes a_jb_i\otimes
ab_j-\sum\limits_{ij}(a_ja_i)b\otimes b_j\otimes
ab_i-\sum\limits_{ij}(b_ib_j)(ab)\otimes a_j\otimes a_i+
$$

$$
\sum\limits_{ij}a_ja_i\otimes b_j\otimes (b_ia)b+
 \sum\limits_{ij}a_j\otimes
a_ib_j\otimes(ab_i)b-\sum\limits_{ij}b_ib_j\otimes(a_ib)a\otimes
a_j-\sum\limits_{ij}a_i\otimes(a_jb)a\otimes b_ib_j-
$$

$$-\sum\limits_{ij}a_ia\otimes b(a_jb_i)\otimes b_j
-\sum\limits_{ij}(a_ja_i)a\otimes bb_i\otimes
b_j+\sum\limits_{ij}a_ja\otimes bb_i\otimes a_ib_j.
$$

The last equality can be rewritten in the form

$$(C_M(r)(b\otimes 1\otimes 1))(a\otimes 1\otimes 1)-C_M(r)(ab\otimes
1\otimes 1)-(C_M(r)(1\otimes 1\otimes a))(1\otimes 1\otimes b)=
$$
$$=C_M(r)(b\otimes 1\otimes a)-C_M(r)(a\otimes b\otimes 1).
$$
This proves the necessary condition.

Let us prove the sufficient condition. Let $r$ satisfies \eqref{UM}.
From the proof of the necessary condition it is clear that
$\Delta_r$  satisfies the second condition of theorem 2. So we need
to prove that $M^*$ is a Malcev algebra and $\Delta_r$ satisfies the
first condition of theorem 1.

In order to prove the first condition we need to check that the
following equality holds:
$$
\sum\limits_i( a_i((ab)c) \otimes
b_i-a_i\otimes((ab)c)b_i)+\sum\limits_i (a_i(bc) \otimes
b_ia-a_i\otimes ((bc)b_i)a)+\sum\limits_i (b(a_i(ac)) \otimes
b_i-ba_i\otimes(ac) b_i)=
$$
$$
=\sum\limits_i ((a_ia)(bc) \otimes b_i-a_i(bc)\otimes
ab_i)-\sum\limits_i (a_ic \otimes (b_ia)b+a_i\otimes ((cb_i)a)b)+
\sum\limits_i ((a_ia)c \otimes b_ib-a_ic\otimes (ab_i)b)+
$$
$$
+\sum\limits_i ((a(a_ib))c \otimes b_i-(aa_i)c\otimes
bb_i)+\sum\limits_i (a_ia \otimes (b_ib)c-a_i\otimes
((ab_i)b)c)+\sum\limits_i (a(a_ib) \otimes b_ic-aa_i\otimes
(bb_i)c)+
$$
$$
\sum\limits_i (a(b(a_ic)) \otimes b_i-a(ba_i)\otimes
cb_i)-\sum\limits_i (a_ib \otimes b_i(ac)+a_i\otimes (bb_i)(ac)).
$$
It holds due to \eqref{mal}.

Let us prove that the dual algebra $M^*$ is a Malcev algebra.

For any $a\in M$ and $f,g,h,t\in M^*$ we have
$$
1. \langle ((fg)h)t,a\rangle =\sum\limits_{ijk} \langle f,
a_k(a_j(a_ia))\rangle\langle g,b_k\rangle\langle h,b_j\rangle\langle
t,b_i\rangle-\sum\limits_{ijk}\langle f,a_k \rangle\langle
g,(a_j(a_ia))b_k \rangle\langle h, b_j \rangle\langle t, b_i
\rangle-
$$
$$
-\sum\limits_{ijk}\langle f, a_ka_j \rangle\langle g,b_k
\rangle\langle h,(a_ia)b_j \rangle\langle t,b_i
\rangle+\sum\limits_{ijk}\langle f,a_k \rangle\langle g,a_jb_k
\rangle\langle h,(a_ia)b_j \rangle\langle t,b_i \rangle-
$$
$$
-\sum\limits_{ijk}\langle f,a_k(a_ja_i) \rangle\langle g,b_k
\rangle\langle h,b_j \rangle\langle t,ab_i
\rangle+\sum\limits_{ijk}\langle f,a_k \rangle\langle g,(a_ja_i)b_k
\rangle\langle h,b_j \rangle\langle t,ab_i \rangle+
$$
$$
+\sum\limits_{ijk}\langle f,a_ka_j \rangle\langle g,b_k
\rangle\langle h,a_ib_j \rangle\langle t,ab_i
\rangle-\sum\limits_{ijk}\langle f,a_k \rangle\langle g,a_jb_k
\rangle\langle h,a_ib_j \rangle\langle t,ab_i \rangle.
$$

$$
2. \langle ((gh)t)f,a\rangle =\sum\limits_{ijk} \langle
f,b_i\rangle\langle g, a_k(a_j(a_ia))\rangle\langle
h,b_k\rangle\langle t,b_j\rangle-\sum\limits_{ijk}\langle f, b_i
\rangle\langle g,a_k \rangle\langle h,(a_j(a_ia))b_k \rangle\langle
t, b_j \rangle-
$$
$$
-\sum\limits_{ijk}\langle f,b_i \rangle\langle g, a_ka_j
\rangle\langle h,b_k \rangle\langle t,(a_ia)b_j
\rangle+\sum\limits_{ijk}\langle f,b_i \rangle\langle g,a_k
\rangle\langle h,a_jb_k \rangle\langle t,(a_ia)b_j \rangle-
$$
$$
-\sum\limits_{ijk}\langle f,ab_i \rangle\langle g,a_k(a_ja_i)
\rangle\langle h,b_k \rangle\langle t,b_j
\rangle+\sum\limits_{ijk}\langle f,ab_i \rangle\langle g,a_k
\rangle\langle h,(a_ja_i)b_k \rangle\langle t,b_j \rangle+
$$
$$
+\sum\limits_{ijk}\langle f,ab_i \rangle\langle g,a_ka_j
\rangle\langle h,b_k \rangle\langle t,a_ib_j
\rangle-\sum\limits_{ijk}\langle f,ab_i \rangle\langle g,a_k
\rangle\langle h,a_jb_k \rangle\langle t,a_ib_j \rangle.
$$

$$
3. \langle ((ht)f)g,a\rangle =\sum\limits_{ijk} \langle
f,b_j\rangle\langle g,b_i\rangle\langle h,
a_k(a_j(a_ia))\rangle\langle t,b_k\rangle-\sum\limits_{ijk}\langle
f, b_j \rangle\langle g, b_i \rangle\langle h,a_k \rangle\langle
t,(a_j(a_ia))b_k \rangle-
$$
$$
-\sum\limits_{ijk}\langle f,(a_ia)b_j \rangle\langle g,b_i
\rangle\langle h, a_ka_j \rangle\langle t,b_k
\rangle+\sum\limits_{ijk}\langle f,(a_ia)b_j \rangle\langle g,b_i
\rangle\langle h,a_k \rangle\langle t,a_jb_k \rangle-
$$
$$
-\sum\limits_{ijk}\langle f,b_j \rangle\langle g,ab_i \rangle\langle
h,a_k(a_ja_i) \rangle\langle t,b_k \rangle+\sum\limits_{ijk}\langle
f,b_j \rangle\langle g,ab_i \rangle\langle h,a_k \rangle\langle
t,(a_ja_i)b_k \rangle+$$
$$
+\sum\limits_{ijk}\langle f,a_ib_j \rangle\langle g,ab_i
\rangle\langle h,a_ka_j \rangle\langle t,b_k
\rangle-\sum\limits_{ijk}\langle f,a_ib_j \rangle\langle g,ab_i
\rangle\langle h,a_k \rangle\langle t,a_jb_k \rangle.
$$

$$
4. \langle ((tf)g)h,a\rangle =\sum\limits_{ijk} \langle
f,b_k\rangle\langle g,b_j\rangle\langle h,b_i\rangle\langle t,
a_k(a_j(a_ia))\rangle-\sum\limits_{ijk}\langle f,(a_j(a_ia))b_k
\rangle\langle g, b_j \rangle\langle h, b_i \rangle\langle t,a_k
\rangle-
$$
$$
-\sum\limits_{ijk}\langle f,b_k \rangle\langle g,(a_ia)b_j
\rangle\langle h,b_i \rangle\langle t, a_ka_j
\rangle+\sum\limits_{ijk}\langle f,a_jb_k \rangle\langle g,(a_ia)b_j
\rangle\langle h,b_i \rangle\langle t,a_k \rangle-
$$
$$
-\sum\limits_{ijk}\langle f,b_k \rangle\langle g,b_j \rangle\langle
h,ab_i \rangle\langle t,a_k(a_ja_i) \rangle+\sum\limits_{ijk}\langle
f,(a_ja_i)b_k \rangle\langle g,b_j \rangle\langle h,ab_i
\rangle\langle t,a_k \rangle+
$$
$$
+\sum\limits_{ijk}\langle f,b_k \rangle\langle g,a_ib_j
\rangle\langle h,ab_i \rangle\langle t,a_ka_j
\rangle-\sum\limits_{ijk}\langle f,a_jb_k \rangle\langle g,a_ib_j
\rangle\langle h,ab_i \rangle\langle t,a_k \rangle.
$$

$$
5. \langle (fh)(gt),a\rangle= \sum\limits_{ijk}\langle f,
a_j(a_ia)\rangle\langle g,a_kb_i\rangle\langle h,b_j\rangle \langle
t,b_k\rangle-\sum\limits_{ijk}\langle f, a_j(a_ia)\rangle\langle g,
a_k\rangle\langle h,b_j \rangle \langle t,b_ib_k \rangle-
$$
$$
-\sum\limits_{ijk}\langle f,a_j \rangle\langle g,
a_kb_i\rangle\langle h, (a_ia)a_j\rangle \langle t,b_k
\rangle+\sum\limits_{ijk}\langle f,a_j \rangle\langle
g,a_k\rangle\langle h,(a_ia)a_j \rangle \langle t,b_ib_k \rangle-
$$
$$-\sum\limits_{ijk}\langle f,a_ja_i
\rangle\langle g,a_k(ab_i)\rangle\langle h,b_j \rangle \langle t,b_k
\rangle+\sum\limits_{ijk}\langle f,a_ja_i \rangle\langle
g,a_k\rangle\langle h,b_j \rangle \langle t,(ab_i)b_k \rangle+
$$
$$
\sum\limits_{ijk}\langle f, a_j \rangle\langle
g,a_k(ab_i)\rangle\langle h,a_ib_j \rangle \langle t,b_k
\rangle-\sum\limits_{ijk}\langle f,a_j \rangle\langle
g,a_k\rangle\langle h,a_ib_j \rangle \langle t,(ab_i)b_k \rangle.
$$

In what follows we will use the following notation: an expression
$x\otimes y\otimes z\otimes s^{i.j}$, $i=1,\dots,5$, $j=1,\dots,8$,
means that $\langle f,x\rangle\langle g,y\rangle\langle
h,z\rangle\langle t,s\rangle$ is equal to the $j$th summand of the
$i$th equality (if $i<5$) or to the negative $j$th summand of the
5th equality (if $i=5$).

Denote by
$$S(a,b)=(C_M(r)(1\otimes b\otimes 1))(1\otimes a\otimes
1)-C_M(r)(ab\otimes 1\otimes 1)-(C_M(r)(1\otimes 1\otimes
a))(1\otimes 1\otimes b)-$$ $$ -C_M(r)(b\otimes 1\otimes
a)+C_M(r)(a\otimes b\otimes 1).$$ Let
 $Q=\sum\limits_i(\xi^2(S(a,a_i)))\otimes b_i=0$. Then by (2) $S(a,b)=0$
for all $a,b\in M$. Consequently, $Q=0$.
Thus (here, expressions
 $p^{(i)}, \ q^{(i)}$ mean, that the elements $p$ and $q$
are equal)
$$
Q= -\sum\limits_{ijk}((a_ja_k)a_i)a\otimes b_k\otimes b_j\otimes
b_i-\sum\limits_{ijk}(b_ka_i)a\otimes b_j\otimes a_ja_k\otimes
b_i^{(1)}- \sum\limits_{ijk}(b_ja_i)a\otimes a_ja_k\otimes
b_k\otimes b_i^{(2)}-
$$
$$
-\sum\limits_{ijk}a_ja_k\otimes b_k \otimes b_j(a_ia)\otimes
b_i^{1.3}-\sum\limits_{ijk} b_k\otimes b_j\otimes
(a_ja_k)(a_ia)\otimes b_i-\sum\limits_{ijk} b_j\otimes a_ja_k\otimes
b_k(a_ia)\otimes b_i^{1.4}+
$$
$$
+\sum\limits_{ijk}a_ja_k\otimes (b_ka)a_i\otimes b_j\otimes
b_i^{5.5}+\sum\limits_{ijk}b_k\otimes (b_ja)a_i\otimes a_ja_k\otimes
b_i^{5.7}+\sum\limits_{ijk}b_j\otimes ((a_ja_k)a)a_i\otimes
b_k\otimes b_i-
$$
$$
-\sum\limits_{ijk}(a_ja_k)a_i\otimes b_k\otimes b_ja\otimes
b_i^{4.6}-\sum\limits_{ijk}b_ka_i\otimes b_j\otimes (a_ja_k)a\otimes
b_i^{(3)}-\sum\limits_{ijk}b_ja_i\otimes a_ja_k\otimes b_ka\otimes
b_i^{4.8}+
$$
$$
+\sum\limits_{ijk}a_ja_k\otimes b_ka\otimes b_ja_i\otimes
b_i^{3.7}+\sum\limits_{ijk}b_k\otimes b_ja\otimes (a_ja_k)a_i\otimes
b_i^{3.5}+\sum\limits_{ijk}b_j\otimes (a_ja_k)a\otimes b_ka_i\otimes
b_i^{(4)}=0.
$$

Let $\varsigma$ be the linear mapping of $M\otimes M\otimes M\otimes
M$ defined by $\varsigma(x\otimes y\otimes z\otimes t)=y\otimes
z\otimes t\otimes x$. Then

$$
\varsigma(Q)=-\sum\limits_{ijk}b_i\otimes ((a_ja_k)a_i)a\otimes
b_k\otimes b_j-\sum\limits_{ijk}b_i\otimes (b_ka_i)a\otimes
b_j\otimes a_ja_k^{(5)}-\sum\limits_{ijk}b_i\otimes (b_ja_i)a\otimes
a_ja_k\otimes b_k^{(4)}-
$$
$$
-\sum\limits_{ijk}b_i\otimes a_ja_k\otimes b_k\otimes
b_j(a_ia)^{2.3}-\sum\limits_{ijk}b_i\otimes b_k\otimes b_j\otimes
(a_ja_k)(a_ia)-\sum\limits_{ijk}b_i\otimes b_j\otimes a_ja_k\otimes
b_k(a_ia)^{2.4}+
$$
$$
+\sum\limits_{ijk}b_i\otimes a_ja_k\otimes (b_ka)a_i\otimes
b_j^{5.3}+\sum\limits_{ijk}b_i\otimes b_k\otimes (b_ja)a_i\otimes
a_ja_k^{5.4}+\sum\limits_{ijk}b_i\otimes b_j\otimes
((a_ja_k)a)a_i\otimes b_k-
$$
$$
-\sum\limits_{ijk}b_i\otimes (a_ja_k)a_i\otimes b_k\otimes
b_ja^{1.6}-\sum\limits_{ijk}b_i\otimes b_ka_i\otimes b_j\otimes
(a_ja_k)a^{(6)}-\sum\limits_{ijk}b_i\otimes b_ja_i\otimes
a_ja_k\otimes b_ka^{1.8}+
$$
$$
+\sum\limits_{ijk}b_i\otimes a_ja_k\otimes b_ka\otimes
b_ja_i^{4.7}+\sum\limits_{ijk}b_i\otimes b_k\otimes b_ja\otimes
(a_ja_k)a_i^{4.5}+\sum\limits_{ijk}b_i\otimes b_j\otimes
(a_ja_k)a\otimes b_ka_i^{(7)}=0.
$$

$$
\varsigma^2(Q)=-\sum\limits_{ijk} b_j\otimes b_i\otimes
((a_ja_k)a_i)a\otimes b_k-\sum\limits_{ijk}a_ja_k\otimes b_i\otimes
(b_ka_i)a\otimes b_j^{(3)}-\sum\limits_{ijk}b_k\otimes b_i\otimes
(b_ja_i)a\otimes a_ja_k^{(7)}-
$$
$$
-\sum\limits_{ijk}b_j(a_ia)\otimes b_i\otimes a_ja_k\otimes
b_k^{3.3}-\sum\limits_{ijk}(a_ja_k)(a_ia)\otimes b_i\otimes
b_k\otimes b_j-\sum\limits_{ijk}b_k(a_ia)\otimes b_i\otimes
b_j\otimes a_ja_k^{3.4}+
$$
$$
+\sum\limits_{ijk}b_j\otimes b_i\otimes a_ja_k\otimes
(b_ka)a_i^{5.8}+\sum\limits_{ijk}a_ja_k\otimes b_i\otimes b_k\otimes
(b_ja)a_i^{5.6}+\sum\limits_{ijk}b_k\otimes b_i\otimes b_j\otimes
((a_ja_k)a)a_i-
$$
$$
-\sum\limits_{ijk}b_ja\otimes b_i\otimes (a_ja_k)a_i\otimes
b_k^{2.6}-\sum\limits_{ijk}(a_ja_k)a\otimes b_i\otimes b_ka_i\otimes
b_j^{(1)}-\sum\limits_{ijk}b_ka\otimes b_i\otimes b_ja_i\otimes
a_ja_k^{2.8}+
$$
$$
\sum\limits_{ijk}b_ja_i\otimes b_i\otimes a_ja_k\otimes
b_ka^{1.7}+\sum\limits_{ijk}(a_ja_k)a_i\otimes b_i\otimes b_k\otimes
b_ja^{1.5}+\sum\limits_{ijk}b_ka_i\otimes b_i\otimes b_j\otimes
(a_ja_k)a^{(8)}=0.
$$

$$
\varsigma^3(Q)=-\sum\limits_{ijk} b_k\otimes b_j\otimes b_i\otimes
((a_ja_k)a_i)a -\sum\limits_{ijk}b_j\otimes a_ja_k\otimes b_i\otimes
(b_ka_i)a^{(6)}-\sum\limits_{ijk}a_ja_k\otimes b_k\otimes b_i\otimes
(b_ja_i)a^{(8)}-
$$
$$
-\sum\limits_{ijk}b_k\otimes b_j(a_ia)\otimes b_i\otimes
a_ja_k^{4.3}-\sum\limits_{ijk}b_j\otimes (a_ja_k)(a_ia)\otimes
b_i\otimes b_k-\sum\limits_{ijk}a_ja_k\otimes b_k(a_ia)\otimes
b_i\otimes b_j^{4.4}+
$$
$$
+\sum\limits_{ijk}(b_ka)a_i\otimes b_j\otimes b_i\otimes
a_ja_k^{5.2}+\sum\limits_{ijk}(b_ja)a_i\otimes a_ja_k\otimes
b_i\otimes b_k^{5.1}+\sum\limits_{ijk}((a_ja_k)a)a_i\otimes
b_k\otimes b_i\otimes b_j-
$$
$$
-\sum\limits_{ijk}b_k\otimes b_ja\otimes b_i\otimes
(a_ja_k)a_i^{3.6}-\sum\limits_{ijk}b_j\otimes (a_ja_k)a\otimes
b_i\otimes b_ka_i^{(5)}-\sum\limits_{ijk}a_ja_k\otimes b_ka\otimes
b_i\otimes b_ja_i^{3.8}+
$$
$$
\sum\limits_{ijk}b_ka\otimes b_ja_i\otimes b_i\otimes
a_ja_k^{2.7}+\sum\limits_{ijk}b_ja\otimes (a_ja_k)a_i\otimes
b_i\otimes b_k^{2.5}+\sum\limits_{ijk}(a_ja_k)a\otimes b_ka_i\otimes
b_i\otimes b_j^{(2)}=0.
$$

Using \eqref{mal} we get
$$-\sum\limits_{ijk}((a_ja_k)a_i)a\otimes b_k\otimes b_j\otimes
b_i-\sum\limits_{ijk} (a_ja_k)(a_ia)\otimes b_i\otimes b_k\otimes
b_j+\sum\limits_{ijk} ((a_ja_k)a)a_i\otimes b_k\otimes b_i \otimes
b_j=$$
$$
=\sum\limits_{ijk} ((a_ia)a_j)a_k\otimes b_k \otimes b_j\otimes
b_i^{1.1}+\sum\limits_{ijk} ((aa_j)a_k)a_i \otimes b_k\otimes  b_j
\otimes b_i^{4.2}$$

Similarly, one can transform the following elements
$$
\sum\limits_{ijk}b_j\otimes ((a_ja_k)a)a_i\otimes b_k\otimes
b_i-\sum\limits_{ijk}b_i\otimes((a_ja_k)a_i)a\otimes b_k\otimes
b_j-\sum\limits_{ijk} b_j\otimes(a_ja_k)(a_ia)\otimes b_i\otimes
b_k=$$
$$=\sum\limits_{ijk}
b_i\otimes a_k(a_j(a_ia))\otimes b_k\otimes
b_j\otimes^{2.1}-\sum\limits_{ijk} a_k\otimes (a_j(a_ia))b_k \otimes
b_j\otimes b_i^{1.2},$$
$$
-\sum\limits_{ijk}b_k\otimes b_j\otimes(a_ja_k)(a_ia)\otimes
b_i+\sum\limits_{ijk}b_i\otimes b_j\otimes ((a_ja_k)a)a_i\otimes
b_k-\sum\limits_{ijk}b_j\otimes b_i\otimes((a_ja_k)a_i)a\otimes b_k=
$$
$$
=\sum\limits_{ijk} b_j\otimes b_i\otimes a_k(a_j(a_ia))\otimes
b_k^{3.1}-\sum\limits_{ijk} b_i\otimes a_k \otimes
(a_j(a_ia))b_k\otimes b_j^{2.2}\text{ and }$$

$$
-\sum\limits_{ijk}b_i\otimes b_k\otimes
b_j\otimes(a_ja_k)(a_ia)+\sum\limits_{ijk}b_k\otimes b_i\otimes
b_j\otimes((a_ja_k)a)a_i-\sum\limits_{ijk}b_k\otimes b_j\otimes
b_i\otimes((a_ja_k)a_i)a=
$$
$$
=\sum\limits_{ijk} b_k\otimes b_j\otimes b_i\otimes
a_k(a_j(a_ia))^{4.1}-\sum\limits_{ijk} b_j\otimes b_i \otimes
a_k\otimes (a_j(a_ia))b_k^{3.2}.$$

Summing up $(id+\varsigma+\varsigma^2+\varsigma^3)Q$, taking into
account the last equalities, and acting on the sum by  $f\otimes
g\otimes h\otimes t$, one can finally get
$$
\langle ((fg)h)t,a\rangle+\langle ((gh)t)f,a\rangle+\langle
((ht)f)g,a\rangle+\langle ((tf)g)h,a\rangle-\langle
(fh)(gt),a\rangle=0.
$$
Thus, the dual algebra  $M^*$ is a Malcev algebra. \hspace*{\fill}
$\Box$

{\bf Corollary 1.} Let $M$  be a Malcev algebra over a field of
characteristic not  2, an element $r\in (id-\tau)(M\otimes M)$ be a
solution to the classical Yang-Baxter equation on $M$. Then the pair
$(M,\Delta_r)$ is a Malcev bialgebra.

\begin{center}
{\bf \S 3. Structures of a Malcev bialgebra on the non-Lie simple
Malcev algebra - preliminary results.}
\end{center}

Let  $F$ be an algebraically closed field of characteristic different from
2 and 3. Then, up to an isomorphism, there is only one non-Lie
simple Malcev algebra $\mathbb{M}$ over $F$ (\cite{Kuz}). The
dimension of $\mathbb{M}$ is equal to seven, and it is convenient to
consider a base
 $h,x,x',y,y',z,z'$ of $\mathbb{M}$
 with the following
multiplication table:
$$
hx=2x,\ hy=2y,\ hz=2z,$$
$$
hx'=-2x',\ hy'=-2y',\ hz'=-2z',$$
$$
xx'=yy'=zz'=h,
$$
$$
xy=2z',\ yz=2x',\ zx=2y',$$
$$x'y'=-2z,\ y'z'=-2x,\
z'x'=-2y.$$ The remaining products are zero. Such a basis is called a
standard basis.

The algebra $\mathbb{M}$ can be constructed in the following way.
Let $\mathcal{C}$ be the matrix Cayley-Dickson algebra with a
multiplication $x\cdot y$. Then $\mathcal{C}=F_2+vF_2$, where $F_2$
is the algebra of $2\times 2$ matrices. Define a new multiplication
in  $\mathcal{C}$: $xy=\frac{1}{2}(x\cdot y-y\cdot x)$. Then the
vector space $\mathcal{C}$ with this multiplication turns into a
Malcev algebra denoted by $\mathcal{C}^{(-)}$. The set $F\cdot 1$,
where $1$ is the unit of $\mathcal{C}$, is the center of
$\mathcal{C}^{(-)}$, and the quotient algebra
$\mathcal{C}^{(-)}/F\cdot1$ is isomorphic to $\mathbb{M}$.

From the multiplication table one can see that the space $M(4)$ with
 the basis $h,x,y',z$ is a four-dimensional subalgebra in $\mathbb{M}$.
  This algebra first appeared in \cite{sgl}.

{\bf Lemma 3.} Suppose $B$ is a subalgebra in $\mathbb{M}$ of
dimension 4. Then  one can choose a standard basis  of $\mathbb{M}$
in such a way that $B=M(4)$.

\textsc{Proof.} The algebra $\mathcal{C}$ can be represented in the
following way: $\mathcal{C}=F\cdot 1+ \mathbb{M}$. The
multiplication in  $\mathcal{C}$ is given by the formula
(\cite{KuzSh})
$$
a\cdot b=-(a,b)1+ab,
$$
 where  $a,b\in
\mathbb{M}$, $(\cdot ,\cdot )$ is a symmetric non-degenerate associative bilinear
form on $\mathbb{M}$ and $ab$ is the antisymmetric multiplication in
$\mathbb{M}$. It is clear that $B'=B+F\cdot 1$ is a five-dimensional
subalgebra in $\mathcal{C}$, and $B$ is isomorphic to the quotient
algebra $B'^{(-)}/F\cdot 1$. Since $B'$ is a non-associative
subalgebra in $\mathcal{C}$, then  $B$  is non-Lie.  In \cite{Gay}
it was proved that over an algebraically closed field of
characteristic different from 2 any four-dimensional non-Lie subalgebra
of $\mathbb{M}$ is isomorphic to $M(4)$. Thus, we can choose a basis
of  $B$ with the same  multiplication table for $h,x,y',z$. So,
we can suppose that  $h,x,y',z$ is the base of $B$. Let $L$ be a
subspace generated by $x,y',z$. Then $(h,L)=(L,L)=0$ and it can
be assumed that $(h,h)=-\frac{1}{4}$. Then the element
$h+\frac{1}{2}\cdot 1$ is an idempotent, so the space $B_1$ with the
base $h+\frac{1}{2}\cdot 1,x,y',z$ is a subalgebra of $\mathcal{C}$,
and $B$ is isomorphic to the quotient algebra $B_1^{(-)}/F\cdot 1$.
It is easy to see that $B_1=F\cdot(h+\frac{1}{2}\cdot 1)+L$ and $L$ is
a nilpotent ideal of $B_1$. We can assume that
$h+\frac{1}{2}\cdot 1 =e_{11}$
(\cite[Lemma 5]{Gme}) and $B_1$ has one of the following
bases:
$e_{11},e_{12},ve_{11}, ve_{12}$
or
$e_{11},ve_{21},ve_{22},e_{21}$,
where $e_{ij}$ is the matrix unit of $F_2$. Then the images of
$$-4e_{11}, 2ve_{11}, 2ve_{22},2e_{21}, 2e_{12}, 2ve_{12},-2ve_{21}$$
under the canonical homomorphism
$\mathcal{C}^{(-)}\mapsto \mathcal{C}^{(-)}/F\cdot 1$ form a
standard basis  of $\mathbb{M}$. It is clear that in this basis
 $B=M(4)$. \hspace*{\fill} $\Box$

Suppose $(\mathbb{M},\Delta)$ is a Malcev bialgebra.
Consider the Drinfeld double $D(\mathbb{M})$. Obviously,
$D(\mathbb{M})$ is not a simple algebra. Therefore $D(\mathbb{M})$
is either a semisimple Malcev algebra or it possesses a non-zero
radical(=solvable radical) $R$. In the following two sections we
consider each of these cases separately.

\begin{center}
{\bf \S 4. The case of non-zero radical.}
\end{center}

{\bf Lemma 4.} Suppose the radical  $R$ of the Drinfeld double
$D(\mathbb{M})$ is nonzero. Then $R=R^{\perp}$, $R^2=0$ and
$D(\mathbb{M})=\mathbb{M}+R$ (semidirect sum).

\textsc{Proof.} By Lemma 1
 $\dim R=\dim\mathbb{M}$. Since $\mathbb{M}$ is a simple algebra and
 $\mathbb{M}$ is not an ideal of $D(\mathbb{M})$, we have
$\mathbb{M}\cap R=0$ and $D(\mathbb{M})=\mathbb{M}+R$. Then,
$\mathbb{M}\cong D(\mathbb{M})/R$.

Let $R^{\perp}$ be the orthogonal complement of $R$ in
$D(\mathbb{M})$ with respect to $Q$. It is straightforward to see
that $R^{\perp}$ is an ideal of
 $D(\mathbb{M})$, so $dim R^{\perp}=dim \mathbb{M}$ by Lemma 2.
 Therefore,
 $D(\mathbb{M})=\mathbb{M}+R^{\perp}$ and $\mathbb{M}$
is isomorphic to the quotient algebra $D(\mathbb{M})/R^{\perp}$.
Let us consider that  $R\cap R^{\perp}= 0$. Then
$D(\mathbb{M})=R\oplus R^{\perp}$ (the direct sum of algebras) and
$R$ is isomorphic to the quotient algebra $D(\mathbb{M})/R^{\perp}$.
It shows that the algebras $\mathbb{M}$ and $R$ are isomorphic what
contradicts the simplicity of $\mathbb{M}$.

Thus  $R\cap R^{\perp}=R$ and $R=R^{\perp}$. Since
$Q(R^2,D(\mathbb{M}))=Q(R,RD(\mathbb{M}))=0$, we finally get
$R^2=0$. \hspace*{\fill} $\Box$

{\bf Theorem 3.} Let a pair  $(\mathbb{M},\Delta)$ be a Malcev
bialgebra, and let the radical $R$ of the Drinfeld double
$D(\mathbb{M})$ is non-zero. Then there exists an element  $r$ from
$(id-\tau)(\mathbb{M}\otimes \mathbb{M})$
such that
$\Delta=\Delta_r$ and $r$ is a solution to the classical Yang-Baxter
equation on $\mathbb{M}$.

\textsc{Proof.} Consider $D(\mathbb{M})$. By Lemma 4
$D(\mathbb{M})=\mathbb{M}+R$, so for every $f\in \mathbb{M}^*$
there exists an element $a\in \mathbb{M}$ such that  $f=a+u$, where
$u\in R$. Define a mapping $\phi\ : \mathbb{M}^* \rightarrow
\mathbb{M}$ by $\phi(f)=a$. The mapping  $\phi$ is a
well-defined homomorphism of algebras. Since $\mathbb{M}^*$ is the
dual space to $\mathbb{M}$, one can find an element
$r=\sum\limits_ia_i\otimes b_i\in \mathbb{M}\otimes \mathbb{M}$,
such that $\phi(f)=\sum\limits_if(b_i)a_i$.

Take $f,g\in \mathbb{M}^*$. Since  $R^{\perp}=R$, then
$Q(f-\phi(f),g-\phi(g))=0$. Consequently,  $$\sum\limits_i
f(b_i)g(a_i)+f(a_i)g(b_i)=0.$$ Therefore,
$$
\sum_i \langle f\otimes g, \sum\limits_i a_i\otimes b_i+b_i\otimes
a_i\rangle=0.$$
Finally, $\tau(r)=-r$.

Further, since $R^2=0$, we get
$(f-\sum\limits_if(a_i)b_i)(g-\sum\limits_i g(a_i)b_i)=0$. Hence,
$$fg-\sum\limits_i
g(b_i)f\leftharpoondown a_i-\sum\limits_i f(b_i)
a_i\rightharpoondown g=0.$$

Thus, for all $a\in \mathbb{M}$
$$
fg(a)=\sum\limits_i \langle f\otimes g, a_ia\otimes b_i+b_i\otimes
aa_i\rangle.
$$

Therefore, $\Delta(a)=\sum\limits_i(a_ia\otimes b_i+b_i\otimes aa_i)$
for all $a\in \mathbb{M}$. Since  $\tau(r)=-r$, then
$\Delta(a)=[r,a]=\sum\limits a_ia\otimes b_i-a_i\otimes ab_i$. In
other words, $\Delta=\Delta_r$.

Since $\phi$ is a homomorphism, we have
$$
\sum\limits_if(b_i)g(b_j)a_ia_j=\sum\limits_ifg(b_i)a_i=\sum\limits_{ij}f(a_jb_i)g(b_j)a_i-\sum\limits_{ij}f(a_j)g(b_ib_j)a_i.
$$
It follows that $\langle f\otimes g \otimes h,
C_{\mathbb{M}}(r)\rangle =0$ for all $f,g,h\in \mathbb{M}^*$.
Therefore $C_{\mathbb{M}}(r)=0$.  \hspace*{\fill} $\Box$

{\bf Theorem 4.} Let $r$ be an antisymmetric solution to the
classical Yang-Baxter equation on  $M$. Then the pair
$(\mathbb{M},\Delta_r)$ is a Malcev bialgebra. Moreover, the radical
$R$ of the Drinfeld double $D(\mathbb{M})$ is nonzero.

\textsc{Proof.} The first statement is valid due to the Theorem 2.
Consider a mapping $\phi: \mathbb{M}^* \rightarrow \mathbb{M}$
defined by $\phi(f)=\sum\limits_if(b_i)a_i$. It is easy to see that
$\phi$ is a homomorphism of algebras. Consider a set  $S=\{
f-\phi(f)|\ f\in \mathbb{M}^*\}$.
For all  $a\in \mathbb{M},f,g\in
\mathbb{M}^*$ we have
$$
(a+g)(f-\phi(f))=-a\phi(f)+a\leftharpoonup f-g\rightharpoonup
\phi(f)+fg+a\rightharpoondown f-g\leftharpoondown \phi(f).
$$
 Since the action of $\phi$ on $fg+a\rightharpoondown
f-g\leftharpoondown \phi(f)$ leads to $a\phi(f) - a\leftharpoonup f
+ g\rightharpoonup \phi(f)$, the right-hand side belongs to $S$,
i.e., $S$ is an ideal in $D(\mathbb{M})$.
%

For all  $f,g\in M^*$ we have
$$Q(f-\phi(f),g-\phi(g))=-f(\phi(g))-g(\phi(f))=\langle f\otimes g,
\sum\limits_i a_i\otimes b_i+b_i\otimes a_i\rangle=0.$$
Using the
associativity of $Q$ we obtain
$$
Q(a+h,(f-\phi(f))(g-\phi(g)))=Q((a+f)(f-\phi(f)),g-\phi(g))=0
$$
for all $a+h\in D(\mathbb{M})$, so $R^2=0$.
Consequently, $S\subseteq R$ and so $R\neq 0$.
\hspace*{\fill} $\Box$

Thus, in order to describe all Malcev bialgebra structures on
$\mathbb{M}$ one should find all antisymmetric solutions of the
classical Yang-Baxter equation  \eqref{YB}. On the other hand,
these solutions (\cite{BD,ZHsimpl,Gme}) are in one to one
correspondence with pairs  $(B,\omega)$, where
$B$ is a subalgebra of $\mathbb{M}$,
$\omega$ is a non-degenerate skew-symmetric
bilinear form satisfying
$$
\omega(xy,z)+\omega(yz,x)+\omega(zx,y)=0
$$
for all $x,y,z\in B$.
In this case $\omega$ is called a symplectic
form, and the pair $(B,\omega)$ is called a symplectic subalgebra.

{\bf Lemma 5.} Let $(B,\omega)$ be a symplectic subalgebra of
$\mathbb{M}$. Then $B$ is isomorphic to one of the following
subalgebra:

1. The subalgebra with the base $x,y'$.

2. The subalgebra with the base $h,x$.

In these cases every non-degenerate skew-symmetric bilinear form is
symplectic.

3. The subalgebra $M(4)$ with the base  $h,x,y',z$. In this case
non-degenerate skew-symmetric bilinear form is symplectic if and
only if it satisfies
$$
\omega(y',h)=2\omega(x,z).
$$

\textsc{ Proof.} Let a pair $(B,\omega)$ be a symplectic subalgebra
of  $\mathbb{M}$. Since $\omega$ is non-degenerate, the
dimension of  $B$ is even. In \cite{Eld}, it was proved that the
dimension of the maximal subalgebra in  $\mathbb{M}$ is equal to 5. Thus,
we have two options for the dimension of $B$: 2 or 4. Let
$\dim B= 2$. In this case,  $B$ is either abelian or
non-abelian solvable Lie algebra. In the first
case, $B$ is isomorphic to the subalgebra with the base $x,y'$,
in the second case---to the subalgebra with the base $h,x$. Clearly,
in the both cases every non-degenerate skew-symmetric bilinear
form is symplectic.

Let $\dim B=4$. By Lemma 3, $B$ is isomorphic to $M(4)$. Let
$\omega$ be a symplectic form on $B$. From the condition
$$
\omega(xz,h)+\omega(zh,x)+\omega(hx,z)=0
$$
we get that  $\omega(y',h)=2\omega(x,z)$. It is easy to see that the
last equality is enough for a non-degenerate skew-symmetric bilinear
form  $\omega$ to be symplectic. \hspace*{\fill} $\Box$

\begin{center}
{\bf \S 5. The semisimple case.}
\end{center}

Consider that $D(M\mathbb{M})$ is a semisimple Malcev algebra. Then
by Lemma 1 $D(\mathbb{M})=M_1\oplus M_2$ is a direct sum of ideals
where each of $M_k$ ($k=1,2$) is isomorphic to $\mathbb{M}$. It is clear
that $Q(M_1,M_2)=0$ and thus the restriction $Q_k$ of the form $Q$
on $M_k$ is a non-degenerate associative form.

 Since $\mathbb{M}^*\subseteq \mathbb{M}+M_1$, then for all  $f\in
\mathbb{M}^*$ we have $f+a\in M_1,$ where $a\in \mathbb{M}$. Let us
define a mapping $\phi_1: \mathbb{M}^* \rightarrow \mathbb{M}$ by
 $\phi_1(f)=-a$. The mapping $\phi_1$ is a well-defined
homomorphism of algebras. Since $M^*$  is a dual space for  $M$,
there is an element
 $r_1=\sum\limits_i a_i\otimes b_i\in \mathbb{M}\otimes \mathbb{M}$,
 such that  $\phi_1(f)=\sum\limits_{i}f(a_i)b_i$.

Similarly, there is an element $r_2=\sum\limits_i c_i\otimes d_i\in
\mathbb{M}\otimes \mathbb{M}$ such that the mapping
 $\phi_2: \mathbb{M}^* \rightarrow \mathbb{M}$ defined by
 $\phi_2(f)=\sum\limits_if(c_i)d_i$ is a homomorphism of algebras
 with $f-\phi_2(f)\in M_2$.
 Since $\phi_1$ and
 $\phi_2$ are homomorphisms, then $r_1$ and $r_2$ are solutions of \eqref{YB1}.

For all $f,g\in \mathbb{M}^*$ we have $Q(f-\phi_1(f),
g-\phi_2(g))=0$. Hence,
$$\sum\limits_i \langle f\otimes g, a_i\otimes b_i+d_i\otimes
c_i\rangle=0.$$ Consequently,
\begin{equation}\label{usl1}
r_1+\tau(r_2)=0.
\end{equation}

Also we have $(f-\phi_1(f))(g-\phi_2(g))=0$.
 Therefore,
 $$fg-f\leftharpoondown \phi_2(g)-
 \phi_1(f)\rightharpoondown g=0.$$
The last equality means that for all $a\in \mathbb{M}$
 $$
 fg(a)=\sum\limits_i f(g(c_i)d_ia)+g(f(a_i)ab_i).$$
 $$
 \langle f\otimes g, \Delta(a)\rangle=\langle f\otimes g,
\sum\limits_i d_ia\otimes c_i+a_i\otimes ab_i\rangle.
$$
Considering  \eqref{usl1} we finally obtain
$$
 \langle f\otimes g, \Delta(a)\rangle=-\langle f\otimes g,
\sum\limits_i a_ia\otimes b_i-a_i\otimes ab_i\rangle.
$$

In other words, we have proved that $\Delta(a)=-(\sum\limits_i
a_ia\otimes b_i - a_i \otimes ab_i)=-\Delta_{r_1}(a)$. Similarly,
one can prove $\Delta(a)=-\Delta_{r_2}(a)$.

Let  $r=r_1$.  Then $\Delta=-\Delta_r$. Since $\mathbb{M}^*$ is an
anticommutative algebra, then
$Q(fg+gf,a)=0$ for all   $a\in \mathbb{M}$.
Therefore,
$$0=\langle f\otimes g, -[r,a]-\tau([r,a])\rangle = -\langle f\otimes g,
[r+\tau(r),a] \rangle.
$$
Thus $[r+\tau(r),a]=0$ for all $a\in \mathbb{M}$.
Put
$s=\frac{1}{2}(r+\tau(r))$, $n=\frac{1}{2}(r-\tau(r)$. Then $r=s+n$,
$n$ is a skew-symmetric element, and $\Delta_r=\Delta_n$.

{\bf Lemma 6.} Let $K$ be an arbitrary Malcev algebra, $r\in
K\otimes K$. Define $s=\frac{1}{2}(r+\tau(r))$,
$n=\frac{1}{2}(r-\tau(r))$. Suppose $r$ is a solution to the
classical Yang-Baxter equation $C_K(r)=0$ with  $[s,a]=0$ for all
$a\in K$. Then the pair $(K,\Delta_r)$ is a Malcev bialgebra.

{\bf Proof.} Since $[s,a]=0$ for all  $a\in K$, then
$\Delta_r=\Delta_n$. It is enough to prove the pair
$(K,\Delta_n)$ to be a Malcev bialgebra. Since  $\tau(n)=-n$, then by
 Theorem 2 it is enough to prove that $C_K(n)$ satisfies
\eqref{UM}.

Let $s=\sum\limits_i a_i\otimes b_i$, $n=\sum\limits_i p_i\otimes
q_i$.  Substituting the sum $s+n$ for $r$ in $C_K(r)=0$,
 we obtain
$$
\sum\limits_{ij}(a_ia_j\otimes b_i\otimes b_j-a_i\otimes
a_jb_i\otimes b_j+a_i\otimes a_j\otimes b_ib_j)+(a_ip_j\otimes
b_i\otimes q_j-a_i\otimes p_jb_i\otimes q_j+a_i\otimes p_j\otimes
b_iq_j)+
$$
$$
+(p_ia_j\otimes q_i\otimes b_j-p_i\otimes a_jq_i\otimes
b_j+p_i\otimes a_j\otimes q_ib_j)+C_K(n)=0.$$

Consider $\sum\limits a_ip_j\otimes b_i\otimes q_j$. Bearing in mind
that  $\sum\limits_i a_i a\otimes b_i+a_i\otimes b_ib=0$, we obtain
$$
\sum\limits_i a_ip_j\otimes b_i\otimes q_j=-\sum\limits_i a_i\otimes
b_ip_j\otimes q_j=\sum\limits_{ij} a_i\otimes p_jb_i\otimes q_j
$$
for all $b\in K$.
Similarly,
$$
\sum\limits_{ij}p_i\otimes a_j\otimes
q_ib_j=\sum\limits_{ij}p_i\otimes a_jq_i\otimes b_j.$$ and
$$
\sum\limits_{ij}a_ia_j\otimes b_i\otimes
b_j=\sum\limits_{ij}a_i\otimes a_jb_i\otimes b_j.
$$
Therefore,
$$
\sum\limits_{ij}a_i\otimes a_j\otimes b_ib_j+a_i\otimes p_j\otimes
b_iq_j+p_ia_j\otimes q_i\otimes b_j+C_K(n)=0.
$$

Consider $\sum\limits_{ij} a_i\otimes p_j\otimes b_iq_j$. Taking
into account that   $\tau(n)=-n$ and $\sum\limits_i a_i a\otimes
b_i+a_i\otimes b_ib=0$ for all $b\in K$ we have
$$
\sum\limits_{ij} a_i\otimes p_j\otimes
b_iq_j=-\sum\limits_{ij}a_iq_j\otimes p_j\otimes
b_i=\sum\limits_{ij}a_ip_j\otimes q_j\otimes b_i=
$$
$$=-\sum\limits_{ij}p_ja_i\otimes q_j\otimes
b_i=-\sum\limits_{ij}p_ia_j\otimes q_i\otimes b_j.$$

Hence we finally obtain
\begin{equation}\label{MYB}
 \sum\limits_{ij}a_i\otimes a_j\otimes
b_ib_j+C_K(n)=0.
\end{equation}

 It follows from \eqref{MYB} that $C_K(n)=-\sum\limits_ia_i\otimes
a_j\otimes b_ib_j$. Plug the expression obtained into \eqref{UM}. For all
$a,b\in K$ we have
$$
(C_K(n)(1\otimes b\otimes 1))(1\otimes a\otimes 1)=-\sum\limits_{ij}
a_i\otimes (a_jb)a \otimes b_ib_j=-\sum\limits_{ij}a_i\otimes
a_j\otimes b_i((b_ja)b)= $$
$$=\sum\limits_{ij}a_i\otimes a_j\otimes
((b_ja)b)b_i,
$$
$$
C_K(n)(ab\otimes 1\otimes 1)=-\sum\limits_{ij}a_i(ab)\otimes
a_j\otimes b_ib_j=\sum\limits_{ij}a_i\otimes a_j\otimes
(b_i(ab))b_j= $$$$=-\sum\limits_{ij}a_i\otimes a_j\otimes
((ab)b_i)b_j,
$$
$$
(C_K(n)(1\otimes 1\otimes a))(1\otimes 1\otimes
b)=-\sum\limits_{ij}a_i\otimes a_j\otimes ((b_ib_j)a)b,
$$
$$
C_K(n)(b\otimes 1\otimes a)=-\sum\limits_{ij} a_ib\otimes a_j\otimes
(b_ib_j)a=\sum\limits_{ij}a_i\otimes a_j\otimes
((b_ib)b_j)a=$$$$-\sum\limits_{ij}a_i\otimes a_j\otimes
((bb_i)b_j)a,
$$
$$
C_K(n)(a\otimes b\otimes 1)=-\sum\limits_{ij} a_ia\otimes
a_jb\otimes b_ib_j=-\sum\limits_{ij} a_i\otimes a_j\otimes
(b_ia)(b_jb)=$$$$=-\sum\limits_{ij} a_i\otimes a_j\otimes
(ab_i)(bb_j),
$$

Finally, by (\ref{mal}), we obtain
$$
(C_K(r)(1\otimes b\otimes 1))(1\otimes a\otimes 1)-C_K(r)(ab\otimes
1\otimes 1)-(C_K(r)(1\otimes 1\otimes a))(1\otimes 1\otimes b)-
$$
$$
-C_K(r)(b\otimes 1\otimes a)+C_K(r)(a\otimes b\otimes 1)=
$$
$$
=\sum\limits_{ij}a_i\otimes a_j\otimes
(((b_ja)b)b_i+((ab)b_i)b_j+((b_ib_j)a)b+((bb_i)b_j)a-(ab_i)(bb_j))=0.$$
\hspace*{\fill} $\Box$

{\bf Lemma 7.} If
 $$l=h\otimes s_1+x\otimes s_2+x'\otimes
s_3+y\otimes s_4+y' \otimes s_5+ z\otimes s_6+z'\otimes s_7\in
\mathbb{M}\otimes \mathbb{M}$$
is such an element that
$[l,a]=0$
for all $a\in \mathbb{M}$
then
$$l=c(\frac{1}{2}h\otimes h+x\otimes x'+x'\otimes
x+y\otimes y'+y'\otimes y+z\otimes z'+z'\otimes z)$$
for some $c\in F$.

{\bf Proof.}  From the condition $[l,h]=0$ we obtain $s_1h=0$
and $-2x\otimes s_2+x\otimes s_2h=0$.
 Therefore, $s_1=ch$ for some
$c\in F$ and $s_2=\beta_1 x'+\beta_2 y'+\beta_3 z'$.

It follows from the condition $[l,x]=0$ that $2cx\otimes
h+x\otimes s_2x=0$. Hence, bearing in mind the condition obtained
for  $s_2$, we finally get $s_2=2c x'$.
Similarly one can prove that
 $s_3=2cx$, $s_4=2c y'$, $s_5=2cy$,
$s_6=2cz'$, $s_7=2cz$. \hspace*{\fill} $\Box$

Since $[s,a]=0$ for all $a\in \mathbb{M}$, then by Lemma 7
$$s=c(\frac{1}{2}h\otimes h+x\otimes x'+x'\otimes x+y\otimes
y'+y'\otimes y+z\otimes z'+z'\otimes z).$$
 If $c=0$ then $\tau(r)=-r$, and by \eqref{usl1}
$r_1=r_2$ and $\phi_1=\phi_2$. But in this case $f+\phi(f)\in
M_1\cap M_2=0$ for all  $f\in \mathbb{M^*}$.

Hence $c\neq0$, so we may assume that $c=\frac{1}{2}$. Then for
$r$ we have
$$
r=h\otimes
(\frac{1}{4}h+\alpha_{12}x+\alpha_{13}x'+\alpha_{14}y+\alpha_{15}y'+\alpha_{16}z+\alpha_{17}z')+
$$
$$
+x\otimes(-\alpha_{12}h+\alpha_{23}x'+\alpha_{24}y+\alpha_{25}y'+\alpha_{26}z+\alpha_{27}z')+
$$
$$
+x'\otimes(-\alpha_{13}h'+(1-\alpha_{23})x+\alpha_{34}y+\alpha_{35}y'+\alpha_{36}z+\alpha_{37}z')+
$$
$$
+y\otimes(-\alpha_{14}h-\alpha_{24}x-\alpha_{34}x'+\alpha_{45}y'+\alpha_{46}z+\alpha_{47}z')+
$$
$$
+y'\otimes(-\alpha_{15}h-\alpha_{25}x-\alpha_{35}x'+(1-\alpha_{45})y+\alpha_{56}z+\alpha_{57}z')+
$$
$$
+z\otimes(-\alpha_{16}h-\alpha_{26}x-\alpha_{36}x'-\alpha_{46}y-\alpha_{56}y'+\alpha_{27}z')+
$$
$$
+z'\otimes(-\alpha_{17}h-\alpha_{27}x-\alpha_{37}x'-\alpha_{47}y-\alpha_{57}y'+(1-\alpha_{27})z).
$$

Let us consider $\Lambda$ to be the following matrix

\begin{equation}\label{mat}
\Lambda=\left ( \begin{array}{ccccccc}\frac{1}{4} &  \alpha_{12} &  \alpha_{13} &  \alpha_{14} &  \alpha_{15} &  \alpha_{16} &  \alpha_{17} \\
-\alpha_{12} & 0 &  \alpha_{23} &  \alpha_{24} &  \alpha_{25} &  \alpha_{26} &  \alpha_{27}\\
-\alpha_{13} & 1-\alpha_{23} &  0 &  \alpha_{34} & \alpha_{35} &
\alpha_{36} &  \alpha_{37}\\
-\alpha_{14} &  -\alpha_{24} &  -\alpha_{34} & 0 &  \alpha_{45} &
\alpha_{46} &  \alpha_{47}\\
-\alpha_{15} &  -\alpha_{25} &  -\alpha_{35} &  1-\alpha_{45} & 0 &
\alpha_{56} &  \alpha_{57}\\
-\alpha_{16} &  -\alpha_{26} &  -\alpha_{36} &  -\alpha_{46} &
-\alpha_{56} & 0 &   \alpha_{67}\\
-\alpha_{17} &  -\alpha_{27} &  -\alpha_{37} &  -\alpha_{47} &
-\alpha_{57} & 1-\alpha_{67} & 0
\end{array} \right).
\end{equation}

We will use the following notations:  $a_1=h,\ a_2=x,\ a_3=x',\
a_4=y\ a_5=y'\ a_6=z\ a_7=z'$. Then
$a=\{a_1,a_2,a_3,a_4,a_5,a_6,a_7\}$ is a basis of $\mathbb{M}$. Let
$\gamma_{ij}^k$ be the structure constants of $\mathbb{M}$ with
respect to  $a_1,\dots, a_7$.
Put
$\Gamma_k=(\gamma_{ij}^k)_{i,j=1,\dots, n}$. Then
$$\Gamma_1=e_{23}-e_{32}+e_{45}-e_{54}+e_{67}-e_{76},$$
$$\Gamma_2=2e_{12}-2e_{21}-2e_{57}+2e_{75},
\Gamma_{3}=-2e_{13}+2e_{31}+2e_{46}-2e_{64},
\Gamma_{4}=2e_{14}-2e_{41}+2e_{37}-2e_{73},$$
$$\Gamma_{5}=2e_{15}-2e_{51}+2e_{26}-2e_{62},
\Gamma_{6}=2e_{16}-2e_{61}-2e_{35}+2e_{53},\Gamma_{7}=-2e_{17}+2e_{71}-2e_{24}+2e_{42}.$$

We have proved that in order to describe all Malcev bialgebra structures
on $\mathbb{M}$ we should find all $r\in \mathbb{M}\otimes
\mathbb{M}$ such that $r+\tau(r)\neq0$, $[r+\tau(r),a]=0$ for all
$a\in \mathbb{M}$ and $C_{\mathbb{M}}(r)=0$.
To proceed in this direction, we need some
properties of the mapping $\phi_1$.

{\bf Lemma 8.}
 Algebras $\mathbb{M}^*$ and $\mathbb{M}$ are not isomorphic.

{\bf Proof.} Assume the converse. Then $\phi_1$ is an isomorphism
of algebras and $\Lambda$ is a non-degenerate matrix.

From \eqref{nd} we obtain (using GAP)
\begin{equation}\label{eq6}
\sum_l2(\Lambda\Gamma_l)_{kl}+(\Lambda^{\tau}\Gamma_k)_{ll}=0.
\end{equation}
Putting $k=1,2,...,7$ one by one into \eqref{eq6}, we get
\begin{equation}\label{e1}
(2\alpha_{23}-1)+(2\alpha_{45}-1)+(2\alpha_{67}-1)=0.
\end{equation}
\begin{equation}\label{e2}
\alpha_{12}-\alpha_{57}=0.
\end{equation}
\begin{equation}\label{e3}
\alpha_{13}-\alpha_{46}=0.
\end{equation}\begin{equation}\label{e4}
\alpha_{14}+\alpha_{37}=0.
\end{equation}\begin{equation}\label{e5}
\alpha_{15}+\alpha_{26}=0.
\end{equation}\begin{equation}\label{e6}
\alpha_{16}-\alpha_{35}=0.
\end{equation}\begin{equation}\label{e7}
\alpha_{17}-\alpha_{24}=0.
\end{equation}

We want to show that there exists an element $p\neq0$ in $\mathbb{M}$
such that $\Delta_r(p)=[r,p]=0$. In this case, the space $\{f\in
\mathbb{M}^*| f(p)=0\}$ is a proper ideal of  $\mathbb{M}^*$,
which contradicts the simplicity of $\mathbb{M}^*$. We have

$$
[r,h]=-4\alpha_{12}(h\otimes x -x\otimes h)+4\alpha_{13}(h\otimes
x'-x'\otimes h)-4\alpha_{14}(h\otimes y-y\otimes
h)+4\alpha_{15}(h\otimes y'-y'\otimes h)-
$$
$$
-4\alpha_{16}(h\otimes z-z\otimes h)+4\alpha_{17}(h\otimes
z'-z'\otimes h)-4\alpha_{24}(x\otimes y-y\otimes
x)-4\alpha_{26}(x\otimes z-z\otimes x)+
$$
$$
+4\alpha_{35}(x'\otimes y'-y'\otimes x')+4\alpha_{37}(x'\otimes
z'-z'\otimes x')-4\alpha_{46}(y\otimes z-z\otimes
x)+4\alpha_{57}(y'\otimes z'-z'\otimes y').
$$

$$
[r,x]=4\alpha_{13}(x\otimes x'-x'\otimes x)-4\alpha_{14}(h\otimes
z'-z'\otimes h-x\otimes y+y\otimes x)+4\alpha_{15}(x\otimes
y'-y'\otimes x)+
$$
$$
+4\alpha_{16}(h\otimes y'-y'\otimes h+x\otimes z-z\otimes
x)+4\alpha_{17}(x\otimes z'-z'\otimes x)+(2\alpha_{23}-1)(h\otimes
x-x\otimes h)-
$$
$$
-4\alpha_{24}(x\otimes z'-z'\otimes x)+4\alpha_{26}(x\otimes
y'-y'\otimes x)-2\alpha_{34}(h\otimes y-y\otimes h+2x'\otimes
z'-2z'\otimes x')-$$$$-2\alpha_{35}(h\otimes y'-y'\otimes h)+
-2\alpha_{36}(h\otimes z-z\otimes h-2x'\otimes y'+2y'\otimes
x')-2\alpha_{37}(h\otimes z'-z'\otimes
h)+$$$$+2(2\alpha_{45}-1)(y'\otimes z'-z'\otimes y')+
4\alpha_{46}(z\otimes z'-z'\otimes z+y\otimes y'-y'\otimes
y)+2(2\alpha_{67}-1)(y'\otimes z'-z'\otimes y').
$$

$$
[r,x']=4\alpha_{12}(x\otimes x'-x'\otimes x)-4\alpha_{14}(x'\otimes
y-y\otimes x')+4\alpha_{15}(h\otimes z-z\otimes h -x'\otimes
y'+y'\otimes x')-
$$
$$
-4\alpha_{16}(x'\otimes z-z\otimes x')-4\alpha_{17}(h\otimes
y-y\otimes h+x'\otimes z'-z'\otimes x')+(2\alpha_{23}-1)(h\otimes
x'-x'\otimes h)+
$$
$$
+4\alpha_{24}(h\otimes y-y\otimes h)+2\alpha_{25}(h\otimes
y'-y'\otimes h+2x\otimes z-2z\otimes x)+4\alpha_{26}(h\otimes
z-z\otimes h)+
$$
$$
2\alpha_{27}(h\otimes z'-z'\otimes h+2x\otimes y-2y\otimes x)
+4\alpha_{35}(x'\otimes z-z\otimes x')-4\alpha_{37}(x'\otimes
y-y\otimes x')+
$$
$$
+2(2\alpha_{45}-1)(y\otimes z-z\otimes y)+4\alpha_{57}(z\otimes
z'-z'\otimes z+y\otimes y'-y'\otimes y)+2(2\alpha_{67}-1)(y\otimes
z-z\otimes y).
$$

$$
[r,y]=4\alpha_{12}(h\otimes z'-z'\otimes h-x\otimes y+y\otimes
x)-4\alpha_{13}(x'\otimes y-y\otimes x')+4\alpha_{15}(y\otimes
y'-y'\otimes y)-
$$
$$
-4\alpha_{16}(h\otimes x'-x'\otimes h-y\otimes z+z\otimes
y)+\alpha_{17}(y\otimes z'-z'\otimes y)-2(2\alpha_{23}-1)(x'\otimes
z'-z'\otimes x')-
$$
$$-4\alpha_{24}(y\otimes z'-z'\otimes y)+2\alpha_{25}(h\otimes x-x\otimes
h-2y'\otimes z'+2z'\otimes y')-$$$$-4\alpha_{26}(z\otimes
z'-z'\otimes z+x\otimes x'-x'\otimes x)+2\alpha_{35}(h\otimes
x'-x'\otimes h)+(2\alpha_{45}+1)(h\otimes y-y\otimes
h)+$$$$+4\alpha_{46}(x'\otimes y-y\otimes x')
-2\alpha_{56}((h\otimes z-z\otimes h-2x'\otimes y'+2y'\otimes
x')-2\alpha_{57}(h\otimes z'-z'\otimes
h)-$$$$-2(2\alpha_{67}-1)(x'\otimes z'-z'\otimes x').
$$

$$
[r,y']=4\alpha_{12}(x\otimes y'-y'\otimes x)-4\alpha_{13}(h\otimes
z-z\otimes h-x'\otimes y'+y'\otimes x')+4\alpha_{14}(y\otimes
y'-y'\otimes y)-
$$
$$
-4\alpha_{16}(y'\otimes z-z\otimes y')+4\alpha_{17}(h\otimes
x-x\otimes h+y'\otimes z'-z'\otimes y')-2(2\alpha_{23}-1)(x\otimes
z-z\otimes x)-
$$
$$
-2\alpha_{24}(h\otimes x-x\otimes h)-4\alpha_{34}(h\otimes
x'-x'\otimes h-y\otimes z+z\otimes y)+4\alpha_{35}(y'\otimes
z-z\otimes y')-
$$
$$
-4\alpha_{37}(x\otimes x'-x'\otimes x+z\otimes z'-z'\otimes
z)+(2\alpha_{45}-1)(h\otimes y'-y'\otimes h)+2\alpha_{46}(h\otimes
z-z\otimes h)+
$$
$$
+2\alpha_{47}(h\otimes z-z\otimes h-x\otimes y+y\otimes
x)-4\alpha_{57}(x\otimes y'-y'\otimes x)-(2\alpha_{67}-1)(x\otimes
z-z\otimes x).
$$

$$
[r,z]=-4\alpha_{12}(h\otimes y'-y'\otimes x+x\otimes z-z\otimes
x)-4\alpha_{13}(x'\otimes z-z\otimes x')+
$$
$$
+4\alpha_{14}(h\otimes x'-x'\otimes h-y\otimes z+z\otimes
y)-4\alpha_{15}(y'\otimes z-z\otimes y')+4\alpha_{17}(z\otimes
z'-z'\otimes z)+
$$
$$
+2(2\alpha_{23}-1)(x'\otimes y'-y'\otimes x') +4\alpha_{24}(x\otimes
x'-x'\otimes x+y\otimes y'-y'\otimes y)-4\alpha_{26}(y'\otimes
z-z\otimes y')+
$$
$$
+2\alpha_{27}(h\otimes x-x\otimes h-2y'\otimes z'+2z'\otimes
y')+2\alpha_{37}(h\otimes x'-x'\otimes
h)+2(2\alpha_{45}-1)(x'\otimes y'-y'\otimes x')+
$$
$$
+4\alpha_{46}(x'\otimes z-z\otimes x')+2\alpha_{47}(h\otimes
y-y\otimes h+2x'\otimes z'-2z'\otimes x')+2\alpha_{57}(h\otimes
y'-y'\otimes h)+$$$$+(2\alpha_{67}-1)(h\otimes z-z\otimes h).
$$

$$
[r,z']=4\alpha_{12}(x\otimes z'-z'\otimes x)+4\alpha_{13}(h\otimes
y-y\otimes h+x'\otimes z'-z'\otimes x')+4\alpha_{14}(y\otimes
z'-z'\otimes y)-
$$
$$
-4\alpha_{15}(h\otimes x-x\otimes h-y'\otimes z'+z'\otimes
y')+4\alpha_{16}(z\otimes z'-z'\otimes z)+2(2\alpha_{23}-1)(x\otimes
y-y\otimes x)-
$$
$$
-2\alpha_{26}(h\otimes x-x\otimes h)+4\alpha_{35}(x\otimes
x'-x'\otimes x+y\otimes y'-y'\otimes y)-2\alpha_{36}(h\otimes
x'-x'\otimes h-2y\otimes z+2z\otimes y)
$$
$$
+4\alpha_{37}(y\otimes z'-z'\otimes y)+2(2\alpha_{45}-1)(x\otimes
y-y\otimes x)-2\alpha_{46}(h\otimes y-y\otimes h)-
$$
$$
-2\alpha_{56}(h\otimes y'-y'\otimes h+2x\otimes z-2z\otimes x)
-4\alpha_{57}(x\otimes z'-z'\otimes x)+(2\alpha_{67}-1)(h\otimes
z'-z'\otimes h).
$$

Take
$p=\alpha_{12}x-\alpha_{13}x'+\alpha_{14}y-\alpha_{15}y'+\alpha_{16}z-\alpha_{17}z'$,
then
$$
[r,\alpha_{12}x-\alpha_{13}x'+\alpha_{14}y-\alpha_{15}y'+\alpha_{16}z-\alpha_{17}z']=
$$
$$
=(\alpha_{12}(2\alpha_{23}-1)+2\alpha_{14}\alpha_{25}+2\alpha_{16}\alpha_{27})(h\otimes
x-x\otimes h-2y'\otimes z'+2z'\otimes y')+$$
$$+(\alpha_{14}(2\alpha_{45}-1)-2\alpha_{12}\alpha_{34}+2\alpha_{16}\alpha_{47})(h\otimes
y-y\otimes h+2x'\otimes z'-2z'\otimes x')+
$$
$$
+(\alpha_{16}(2\alpha_{67}-1)-2\alpha_{12}\alpha_{36}-2\alpha_{14}\alpha_{56})(h\otimes
z-z\otimes h-2x'\otimes y'+2y'\otimes x')+
$$
$$
+(2\alpha_{15}\alpha_{34}+2\alpha_{17}\alpha_{36}-\alpha_{13}(2\alpha_{23}-1))(h\otimes
x'-x'\otimes h-2y\otimes z+2z\otimes x)+
$$
$$
+(2\alpha_{17}\alpha_{56}-2\alpha_{13}\alpha_{25}-\alpha_{15}(2\alpha_{45}-1))(h\otimes
y'-y'\otimes h+2x\otimes z-2z\otimes x)+
$$
$$
+(-2\alpha_{13}\alpha_{27}-2\alpha_{15}\alpha_{47}-\alpha_{17}(2\alpha_{67}-1))(h\otimes
z'-z'\otimes h-2x\otimes y+2y\otimes x).
$$

In \eqref{eq2} put $k=7, s=1, n=5$ to get
$$
0=\alpha_{14}\alpha_{25}+\alpha_{16}\alpha_{27}-\alpha_{12}\alpha_{
45}-\alpha_{12}\alpha_{67}+\alpha_{12}=\alpha_{14}\alpha_{25}+\alpha_{16}\alpha_{27}+\alpha_{12}((\alpha_{45}-\frac{1}{2})+(\alpha_{67}-\frac{1}{2}))=
$$
$$
=\alpha_{14}\alpha_{25}+\alpha_{16}\alpha_{27}+\alpha_{12}(\alpha_{23}-\frac{1}{2}).
$$
Thus, the coefficient of the first summand of $[r,p]$ is equal to zero.
Similarly, putting
  $\{k=7, s=1, n=3\}$, $\{k=3, s=1, n=5\}$, $\{k=6, s=1,
n=4\}$, $\{k=6, s=1, n=2\}$, $\{k=2, s=1, n=4\}$, we obtain
$[r,p]=0$.\hspace*{\fill} $\Box$

By Lemma 2 the space
$$
V=\{a\in \mathbb{M}|\ a+f\in M_2\ \text{for some}\ f\in
\mathbb{M}^*\}
$$
is a subbialgebra of $(\mathbb{M},\Delta)$ and $M_2V^\perp=0$.
Then $V^\perp\subseteq M_1$. By Lemma 8  $V\neq \mathbb{M}$.
In \cite{Eld} it was proved
that any maximal subalgebra in $\mathbb{M}$ is isomorphic to the
algebra $M(5)$ with base $h,x,x',y',z$.

On the other hand, if $\dim V\leq 3$ then  $\dim(\mathbb{M}^*\cap
M_1)\geq 4$. But $Q(\mathbb{M}^*,\mathbb{M}^*)=0$ that is impossible
since $Q_1$ is non-degenerate. So, we have two options: either $\dim V=4$
or $\dim V=5$.

{\bf Lemma 9.}  The dimension of $V$ does not equal 5.

\textsc{Proof.} Suppose that $\dim V=5$. We can assume, that the
elements $h,x,x',y',z$ form the basis of $V$. By Lemma 2 the space
$$V^{\perp}=\{f\in M^*|\ f(a)=0\ \text{for every}\ a\in V\}$$
satisfies $M_2V^{\perp}=0$. So, $V^{\perp}\subseteq M_1$ and,
therefore,  $V^{\perp}=ker \phi_1$. If
$h^*,x^*,x'^*,y^*,y'^*,z^*,z'^*$ is a dual base
of $\mathbb{M}^*$ to the base
$h,x,x',y,y',z,z'$  then the functionals $y^*$ and
$z'^*$ form a base of $V^{\perp}$. In this case,
some entries of the matrix
\eqref{mat} are zero:
$$\alpha_{14}=\alpha_{24}=\alpha_{34}=\alpha_{45}=\alpha_{46}=\alpha_{47}=0,
\alpha_{j7}=0$$ for all $i,j=1,\dots, 6$.
Thus
\begin{equation}\label{mat2}
\Lambda=\left ( \begin{array}{ccccccc}\frac{1}{4} &  \alpha_{12} &  \alpha_{13} & 0 &  \alpha_{15} &  \alpha_{16} &  0 \\
-\alpha_{12} & 0 &  \alpha_{23} &  0 &  \alpha_{25} &  \alpha_{26} &  0\\
-\alpha_{13} & 1-\alpha_{23} &  0 &  0 & \alpha_{35} &
\alpha_{36} &  0\\
0 &  0 &  0 & 0 &  0 & 0 &  0\\
-\alpha_{15} &  -\alpha_{25} &  -\alpha_{35} &  1 & 0 &
\alpha_{56} &  0\\
-\alpha_{16} &  -\alpha_{26} &  -\alpha_{36} &  0 &
-\alpha_{56} & 0 &   1\\
0 &  0 &  0 &  0 & 0 & 0 & 0
\end{array} \right)
\end{equation}
Since $dim V=5$, then the rank of $\Lambda$ is equal to 5.

Now, if we put $\Lambda$ into \eqref{eq2}, we obtain the following
equalities (by means of GAP):

\begin{equation}\label{e8}
4\alpha_{12}\alpha_{13}-\alpha_{23}(\alpha_{23}-1)=0\ \
\{n=3,k=2,s=1\};
\end{equation}

\begin{equation}\label{e9}
4\alpha_{12}\alpha_{15}+2\alpha_{12}\alpha_{26}-\alpha_{25}(\alpha_{23}-1)=0\
\ \{n=5,\ k=2,\ s=1\};
\end{equation}

\begin{equation}\label{e10}
2\alpha_{12}\alpha_{36}+2\alpha_{16}(\alpha_{23}-1)-\alpha_{35}(\alpha_{23}-1)=0\
\ \{n=5,\ k=3,\ s=1\};
\end{equation}

\begin{equation}\label{e11}
-4\alpha_{12}\alpha_{35}+2\alpha_{26}(\alpha_{23}-1)=0\ \ \{n=5,\
k=3,\ s=2\};
\end{equation}

\begin{equation}\label{e12}
2\alpha_{13}\alpha_{25}-2\alpha_{15}\alpha_{23}-\alpha_{23}\alpha_{26}=0\
\ \{n=6,\ k=2,\ s=1\};
\end{equation}

\begin{equation}\label{e13}
-4\alpha_{13}\alpha_{16}+2\alpha_{13}\alpha_{35}-\alpha_{23}\alpha_{36}=0\
\ \{n=6,\ k=3,\ s=1\};
\end{equation}

\begin{equation}\label{e14}
-4\alpha_{13}\alpha_{26}+2\alpha_{23}\alpha_{35}=0\ \ \{n=6,\ k=3,\
s=2\};
\end{equation}

\begin{equation}\label{e15}
-4\alpha_{15}\alpha_{16}+2\alpha_{15}\alpha_{35}-2\alpha_{16}\alpha_{26}-\alpha_{25}\alpha_{36}+\alpha_{26}\alpha_{35}=0\
\ \{n=6,\ k=5,\ s=1\};
\end{equation}

\begin{equation}\label{e16}
-4\alpha_{15}\alpha_{26}+2\alpha_{25}\alpha_{35}-2\alpha_{26}^2=0\ \
\{n=6,\ k=5,\ s=2\};
\end{equation}

\begin{equation}\label{e17}
-4\alpha_{16}\alpha_{35}-2\alpha_{26}\alpha_{36}+2\alpha_{35}^2=0\ \
\{n=6,\ k=5,\ s=3\};
\end{equation}

Other relations coming from another values of $n,k,s$
follow from the given equalities. Consider a matrix
$$
\Lambda_1=\left ( \begin{array}{ccccccc}\frac{1}{4} &  \alpha_{12} &  \alpha_{13} & 0 &  \alpha_{15} &  \alpha_{16} &  0 \\
-\alpha_{12} & 0 &  \alpha_{23} &  0 &  \alpha_{25} &  \alpha_{26} &  0\\
-\alpha_{13} & 1-\alpha_{23} &  0 &  0 & \alpha_{35} & \alpha_{36} &
0
\end{array} \right).
$$
Let us prove that the rank of $\Lambda_1$ does not exceed 2. For
this, consider a matrix
$$
\Lambda_2=\left ( \begin{array}{ccc}\frac{1}{4} &  \alpha_{12} &  \alpha_{13}  \\
-\alpha_{12} & 0 &  \alpha_{23}\\
-\alpha_{13} & 1-\alpha_{23} &  0
\end{array} \right).
$$
Let $$
\Lambda^*_2=\left ( \begin{array}{ccc}-\alpha_{23}(1-\alpha_{13}) &  \alpha_{13}(1-\alpha_{23}) &  \alpha_{12}\alpha_{23}  \\
-\alpha_{13}\alpha_{23} & \alpha_{13}^2 &  -(\frac{1}{4}\alpha_{23}+\alpha_{12}\alpha_{13})\\
-\alpha_{12}(1-\alpha_{23}) &
-(\frac{1}{4}(1-\alpha_{23})+\alpha_{12}\alpha_{13}) & \alpha_{12}^2
\end{array} \right).
$$ Then
$\Lambda_2\Lambda_2^*=\det(\Lambda_2)E$.

Let $V_2$ and $V_3$ be  second and third rows of $\Lambda_2^*$,
respectively. It is easy to see that $V_2$ and $V_3$ can not be
equal to zero simultaneously. By \eqref{e8} $\det(\Lambda_2)=0$, so
$\Lambda_2\Lambda_2^*=0$. In particular, $V_i\Lambda_2=0$, for
$i=2,3$.

Let $U_5$ and $U_6$ be  5th and 6th columns of $\Lambda_1$,
respectively. Let us prove that  $V_iU_j=0$ for all $i=2,3,\ j=5,6$.
We have
$$
V_2U_5=-\alpha_{13}\alpha_{23}\alpha_{15}+\alpha_{13}^2\alpha_{25}-(\frac{1}{4}\alpha_{23}+\alpha_{12}\alpha_{13})\alpha_{35}=
^{\eqref{e12}}-\alpha_{13}\alpha_{23}\alpha_{15}+\alpha_{13}\alpha_{23}\alpha_{15}+$$$$+\frac{1}{2}\alpha_{13}\alpha_{23}\alpha_{26}
-(\frac{1}{4}\alpha_{23}+\alpha_{12}\alpha_{13})\alpha_{35}=^{\eqref{e14}}\alpha_{35}(\frac{1}{4}\alpha_{23}(\alpha_{23}-1)-\alpha_{12}\alpha_{13})=
^{\eqref{e8}}0.
$$

$$
V_2U_6=-\alpha_{13}\alpha_{23}\alpha_{16}+\alpha_{13}^2\alpha_{26}-(\frac{1}{4}\alpha_{23}+\alpha_{12}\alpha_{13})\alpha_{36}=^{\eqref{e14}}
\alpha_{13}(-\alpha_{23}\alpha_{16}+\frac{1}{2}\alpha_{23}\alpha_{35}-\alpha_{12}\alpha_{36})-
$$
$$
-\frac{1}{4}\alpha_{23}\alpha_{36}=^{\eqref{e10}}-\alpha_{13}\alpha_{16}+\frac{1}{2}\alpha_{13}\alpha_{35}-\frac{1}{4}\alpha_{23}\alpha_{36}=^{\eqref{e13}}0
$$

Similarly, $V_3U_5=V_3U_6=0$. It follows that the rows of
$\Lambda_1$ are linearly dependent. Therefore the rank of $\Lambda$
does not exceed 4. \hspace*{\fill} $\Box$

Now we are ready to prove the main theorem in this section.

{\bf Theorem 5.} Let $\mathbb{M}$ be a simple non-Lie Malcev algebra
over an algebraically closed field of characteristic not equal 2, 3.
Then in a standard basis $h,x,x',y,y',z,z'$ of the algebra
$\mathbb{M}$
  the element
$$ r_0=\alpha_{12} (h\otimes x-x\otimes h)+ \alpha_{15}(h\otimes
y'-y'\otimes h)+\alpha_{16}(h\otimes z-z\otimes h) +$$ $$
\alpha_{25}(x\otimes y'-y'\otimes x)- 2\alpha_{15}(x\otimes
z-z\otimes x)+\alpha_{56}(y'\otimes z-z\otimes y')$$  is a solution
of the classical Yang-Baxter equation on $\mathbb{M}$. Moreover, the
element
\begin{equation}\label{us4} r=r_0+\frac{1}{4}h\otimes
h+x\otimes x'+y'\otimes y+z\otimes z'.\end{equation} induces on
$\mathbb{M}$ a structure of a Malcev bialgebra with a semisimple
Drinfeld double.

Conversely, let  $(\mathbb{M},\Delta)$ be a Malcev bialgebra with a
semisimple Drinfeld double. Then $\Delta=-\Delta_r$, and one can
choose a standard basis $h,x,x',y,y',z,z'$ of  $\mathbb{M}$
in such  a way that
$ r$ has the form (\ref{us4}).

{\bf Proof.} By Lemma 5 and (\ref{MYB}) we obtain $r_0$ to be a
solution to the classical Yang-Baxter equation on $\mathbb{M}$.

Let $r=r_0+\frac{1}{4}h\otimes h+x\otimes x'+y'\otimes y+z\otimes
z'$. For $r$ to be a solution to the equation $C_{\mathbb{M}}(r)=0$
it is necessary and sufficient that the equalities
\eqref{e8}--\eqref{e17} hold. Since in our case
$\alpha_{13}=\alpha_{35}=\alpha_{36}=0$ and $\alpha_{23}=1$,
straightforward calculations (using GAP) show that
\eqref{e8}--\eqref{e17} follow from $2\alpha_{15}+\alpha_{26}=0$.
Thus, by Lemma 6, the pair
 $(\mathbb{M},\Delta_r)$ is a Malcev bialgebra. Since
$r\neq -\tau(r)$, then the Drinfeld double $D(M)$ has to be a
semisimple algebra.

Conversely, let $(\mathbb{M},\Delta)$ be a Malcev bialgebra with a
semisimple Drinfeld double. By Lemmas 8 and 9 the dimension of $V$
equals 4. Thence,
 $\ker\phi_1=V^{\perp}$.
Then, by Lemma 3 one can choose  a standard basis $h,x,x',y,y',z,z'$
of $\mathbb{M}$ in such a way that  $V$ has the base $h,x,y',z$. If
$h^*,x^*,x'^*,y^*,y'^*,z^*,z'^*$ is the dual base to the base
$h,x,x',y,y',z,z'$  then the elements $x'^*,y^*,z'^*$ form a base of
$V^{\perp}$. Therefore,
\begin{equation}
\Lambda=\left ( \begin{array}{ccccccc}\frac{1}{4} &  \alpha_{12} &  0 & 0 &  \alpha_{15} &  \alpha_{16} &  0 \\
-\alpha_{12} & 0 &  1 &  0 &  \alpha_{25} &  \alpha_{26} &  0\\
0 & 0 &  0 &  0 & 0 &
0 &  0\\
0 &  0 &  0 & 0 &  0 & 0 &  0\\
-\alpha_{15} &  -\alpha_{25} &  0 &  1 & 0 &
\alpha_{56} &  0\\
-\alpha_{16} &  -\alpha_{26} &  0 &  0 &
-\alpha_{56} & 0 &   1\\
0 &  0 &  0 &  0 & 0 & 0 & 0
\end{array} \right).\nonumber
\end{equation}
From the equality \eqref{eq2} with $n=6,\ k=2,\ s=1$ we obtain
$2\alpha_{15}+\alpha_{26}=0$. Hence $r$ has the form \eqref{us4}.
 \hspace*{\fill} $\Box$

\renewcommand{\refname}{\begin{center} References \end{center}}

\end{document}